\documentclass[11pt]{amsart}

\issueinfo{00}{0}{Month}{Year}
\copyrightinfo{2008}{University of Calgary}
\pagespan{1}{2}
\PII{ISSN 1715-0868}
\include{cdmstdcmds}

  \usepackage{amsmath}
  \usepackage{amssymb}
  \usepackage{cite}
  \usepackage[utf8]{inputenc}
  \usepackage{amsthm}
  \usepackage[hidelinks]{hyperref}
  \usepackage{graphicx}
  \usepackage{tikz-cd}
  \usepackage{color}
  \usepackage[justification=centering]{caption}

  \newcommand{\N}{\mathbb{N}}
  \newcommand{\Z}{\mathbb{Z}}
  \newcommand{\R}{\mathbb{R}}
  \newcommand{\C}{\mathcal{C}}
  
  \newcommand{\tX}{\widetilde{X}}

  \newtheorem{thm}{Theorem}[subsection]
  \newtheorem{lemma}[thm]{Lemma}
  \newtheorem{prop}[thm]{Proposition}

  \theoremstyle{definition}
  \newtheorem{df}[thm]{Definition}
  \newtheorem{ex}[thm]{Example}

  \newtheorem{nota}[thm]{Notation}

  \theoremstyle{remark}

\begin{document}

\title{The Lifting Properties of A-Homotopy Theory}

\author{Rachel Hardeman Morrill}
\address{Department of Mathematics and Statistics\\University of Calgary\\2500 University Dr.~NW\\Calgary, Alberta, Canada T2N 1N4}
\email{rachel.hardeman1@ucalgary.ca}

%

\date{(date1), and in revised form (date1).}

\thanks{}

\begin{abstract}
In classical homotopy theory, two spaces are considered homotopy equivalent if one space can be continuously deformed into the other. This theory, however, does not respect the discrete nature of graphs. For this reason, a discrete homotopy theory that recognizes the difference between the vertices and edges of a graph was invented, called A-homotopy theory \cite{Atkin1, Atkin2, BabsonHomotopy, BarceloFoundations, BarceloPerspectives}. In classical homotopy theory, covering spaces and lifting properties are often used to compute the fundamental group of the circle. In this paper, we develop the lifting properties for A-homotopy theory. Using a covering graph and these lifting properties, we compute the fundamental group of the cycle $C_{5}$, giving an alternate approach to \cite{BarceloFoundations}.
\end{abstract}

\maketitle


\section{Introduction}
A-homotopy theory is a discrete homotopy theory developed to investigate the invariants of graphs in a combinatorial way that respects the stucture of graphs \cite{BarceloPerspectives}. The first iteration of this theory is called Q-analysis, developed by Ron Atkin as a means to study the combinatorial “holes” of simplicial complexes \cite{Atkin1, Atkin2}. Q-analysis is used in several fields, including sociology and biology \cite{BarceloFoundations}. Atkin noticed that Q-analysis was the foundation of a general theory and gave a road map for constructing it \cite{BarceloPerspectives}. Barcelo et al. developed this general theory for simplicial complexes and a related theory specifically for graphs in \cite{BarceloFoundations}. This related theory, named A-homotopy theory in honor of Aktin, was further developed by Babson et al. in \cite{BabsonHomotopy}. It allows us to find areas where there are fewer edges connecting the vertices of a  graph. Since graphs are often used to represent real world networks and systems, these areas with fewer edges can either point to missing information in the network or areas where the network could be made more efficient by adding connections.
    
Throughout algebraic topology, several tools are used to examine topological spaces. Algebraic topology is a century-old theory, and the definitions and techniques are well established in the literature. They can be found in any algebraic topology textbook, including \cite{HatcherAlgebraic, MunkresTopology}. The \textit{classical fundamental group of a space}, also called the first homotopy group, is defined to be the set of continuous maps from the unit interval into the space under an equivalence. These continuous maps send both endpoints of the interval to the base point of the space. We compute the fundamental group of spaces using several different methods, including covering spaces and lifting properties. A \textit{covering space of a space $A$} is a space $\widetilde{A}$ together with a continuous onto map $\rho: \widetilde{A} \to A$ which is locally a homeomorphism. That is, the space $\widetilde{A}$ looks like the space $A$ locally. Given a covering space $\rho: \widetilde{A} \to A$ and a continuous map $f: B \to A$, a \textit{lift} is a map $\widetilde{f}: B \to \widetilde{A}$ which factors $f$ through the map $\rho$. In algebraic topology, the lifting properties tell us under what conditions these lifts exist.

In \cite[Proposition 5.12]{BarceloFoundations}, Barcelo et al. prove that attaching 2-cells to the 3-cycles and 4-cycles of a graph, and computing the classical fundamental group of the space created, is equivalent to computing the A-homotopy fundamental group of the original graph. As a consequence of this, the cycles $C_3$ and $C_4$ are contractible as graphs. This means that they are homotopy equivalent to a single vertex, while all cycles $C_k$ for $k \geq 5$ are not contractible. From a topologist’s perspective this is unusual, because all cycles $C_k$ for $k \geq 3$, including $C_3$ and $C_4$, are homotopy equivalent to the circle as spaces, and thus not contractible in classical homotopy theory. In \cite{BabsonHomotopy}, Babson et al. construct a cubical set for each graph. Every cubical set has a topological space associated to it called the geometric realization. In \cite[Theorem 5.2]{BabsonHomotopy}, Babson et al. prove that computing the classical homotopy groups of the geometric realization associated to the cubical set of a graph is equivalent to computing the A-homotopy groups of the original graphs. Both of these computation methods require the use of computational techniques in classical homotopy theory, outside of graph theory. 

To make A-homotopy theory more self-contained, we develop an approach to computing the A-homotopy fundamental group of a graph in this paper analogous to the covering space method used in classical homotopy theory. Covering spaces and lifting properties are important for understanding the fundamental group in topology. This approach offers a different perspective that does not involve attaching 2-cells to graphs and using classical homotopy theory, but instead allows us to stay in the world of graphs. We will use the existing definition of covering graphs found in \cite{GodsilAlgebraic} and prove the Path Lifting Property (Theorem \ref{Path_Lifting_Property}), which allows us to factor graph homomorphisms from paths into a graph $X$ through a covering graph of $X$. We use the Path Lifting Property to prove the Homotopy Lifting Property (Theorem \ref{Homotopy_Lifting_Property}), which allows us to factor a graph homotopy through a covering graph. Finally, we use both the Path Lifting Property and the Homotopy Lifting Property to prove the Lifting Criterion (Theorem \ref{Lifting_Criterion}), which gives us the conditions for when an arbitrary graph homomorphism can be factored through a covering graph. There are interesting conditions on the Homotopy Lifting Property (Theorem \ref{Homotopy_Lifting_Property}) and the Lifting Criterion (Theorem \ref{Lifting_Criterion}) because 3-cycles and 4-cycles are contractible in A-homotopy. 

To demonstrate that covering graphs are effective computational tools, we compute the A-homotopy fundamental group of the $k$-cycle for $k \geq 5$ using a covering graph and lifting properties in a method analogous to the method used to compute the classical fundamental group of the circle. This method recovers the computation of the fundamental group of $C_k$ for $k \geq 5$ in \cite{BarceloFoundations} as expected. While covering graphs can be used as a computational tool, that is not the reason we are developing this theory. Covering spaces are an important part of classical homotopy theory and are an important special case of Hurewicz fibrations. These fibrations are used to construct a homotopy category in topology. We will further develop the covering graph theory, including the universal covers, for A-homotopy theory in another paper.

This paper is organized as follows. In section 2, we give the basic definitions and theorems of A-homotopy theory currently found in the literature \cite{BabsonHomotopy}. In section 3, we give the definition of covering graphs found in \cite{HellGraphs} and develop the Path Lifting Property (Theorem \ref{Path_Lifting_Property}), the Homotopy Lifting Property (Theorem \ref{Homotopy_Lifting_Property}), and
the Lifting Criterion (Theorem \ref{Lifting_Criterion}) for A-homotopy theory. In section 4, we use a covering graph and these lifting properties to provide an alternate approach to computing the fundamental group of the $k$-cycle for $k \geq 5$.

\subsection{Acknowledgements}
The results from this paper come from my thesis \cite{HardemanMscThesis}. Many
thanks to my masters supervisors Professor Kristine Bauer and Professor Karen
Seyffarth of the University of Calgary for your endless support and kindness. I
would also like to thank Professor H{\'e}l{\`e}ne Barcelo of the Mathematical Sciences Research Institute, Berkeley, California for making the time to meet with me and discuss the development of A-homotopy theory. Thank you to NSERC and PIMS as well for the funding that supported this research.

\section{Background}
In this section, we review the fundamental definitions of A-homotopy theory. We proceed by analogy with the basic definition of classical homotopy theory and discuss the major differences between these two theories. We claim no originality. This theory was developed in \cite{BabsonHomotopy, BarceloFoundations, BarceloPerspectives}. But first, we need to consider some basic definitions and lemmas that are the building blocks of this discrete homotopy theory.

\subsection{Basic Definitions}
In classical homotopy theory, the main objects are topological space, and we examine continuous maps between topological spaces. We frequently use the product of two spaces. We continuously deform maps over the unit interval, and when forming the fundamental group, we map the unit interval into the space. In A-homotopy theory, the main objects are simple graphs with a loop attached to every vertex, and we use graph homomorphisms between these graphs. These graph homomorphisms respect the structure of our graphs. We employ a standard discrete product found in graph theory, called the Cartesian product. In order to better distinguish between vertices and edges in the graphs that we examine, we replace the unit interval with graphs known as paths of length $n$. These paths are denoted $I_n$ with vertices labelled from 0 to a non-negative integer $n$. We also use the path of infinite length, denoted by $I_{\infty}$, with vertices labelled by the integers.

While \cite{BabsonHomotopy, BarceloFoundations, BarceloPerspectives} used simple graphs as their main objects, we have attached a loop to every vertex so that we can use the standard definition of graph homomorphism from graph theory. For a graph $X$, we denote the set of vertices of $X$ by $V(X)$, the set of edges by $E(X)$,
and an edge between the vertices $v$ and $w$ by $vw$. If $vw \in E(X)$, then we say that the vertices $v$ and $w$ are \textit{adjacent}. Some graphs we consider have one selected vertex called a \textit{distinguished vertex}. We denote a graph $X$ with distinguished vertex $x_0$ by $(X, x_0)$, and we call these based graphs. In later sections, we will examine paths in a graph $X$ that start and end at the distinguished vertex of $X$. 

\begin{df}\cite[p. 6]{GodsilAlgebraic}\label{Graph_Homomorphism_Def}
A \textit{graph homomorphism} $f: X \to Y$ is a map of sets $V(X) \to V(Y)$ such that, if $uv \in E(X)$, then $f(u)f(v) \in E(Y)$, that is, the image of a pair of adjacent vertices in $X$ is a pair of adjacent vertices in $Y$.
\end{df}

Since we use graphs with a loop at every vertex, we can still map two adjacent vertices to the same vertex, as every vertex is adjacent to itself. For example, there is a graph homomorphism from any graph $X$ into the graph with a single vertex.

\begin{ex}\label{Graph_Homomorphism_Ex}
In Figure \ref{Graph_Homomorphism_Fig}, the graph on the left is denoted $(C_3, x)$ and the graph on
the right is denoted $(R, a)$. The distinguished vertices are shown in green. For the sake of clarity, all figures in this paper will not display the loops at every vertex. The vertex set maps $f: C_3 \to R$ and $g: R \to C_3$ defined
by
\[ \begin{array}{ccc}
f(x) = a, \;\;\; & & g(a) = x, \\
f(y) = d, \;\;\; & \textrm{and} \;\;\; & g(b) = y, \\
f(z) = c, \;\;\; & & g(c) = z, \\
\;\;\; & & g(d) = y.
\end{array} \]
are examples of graph homomorphisms.
\begin{figure}[ht]
	\begin{picture}(4.5, 1.9)
	\put(0, 0){\includegraphics[scale=0.5]{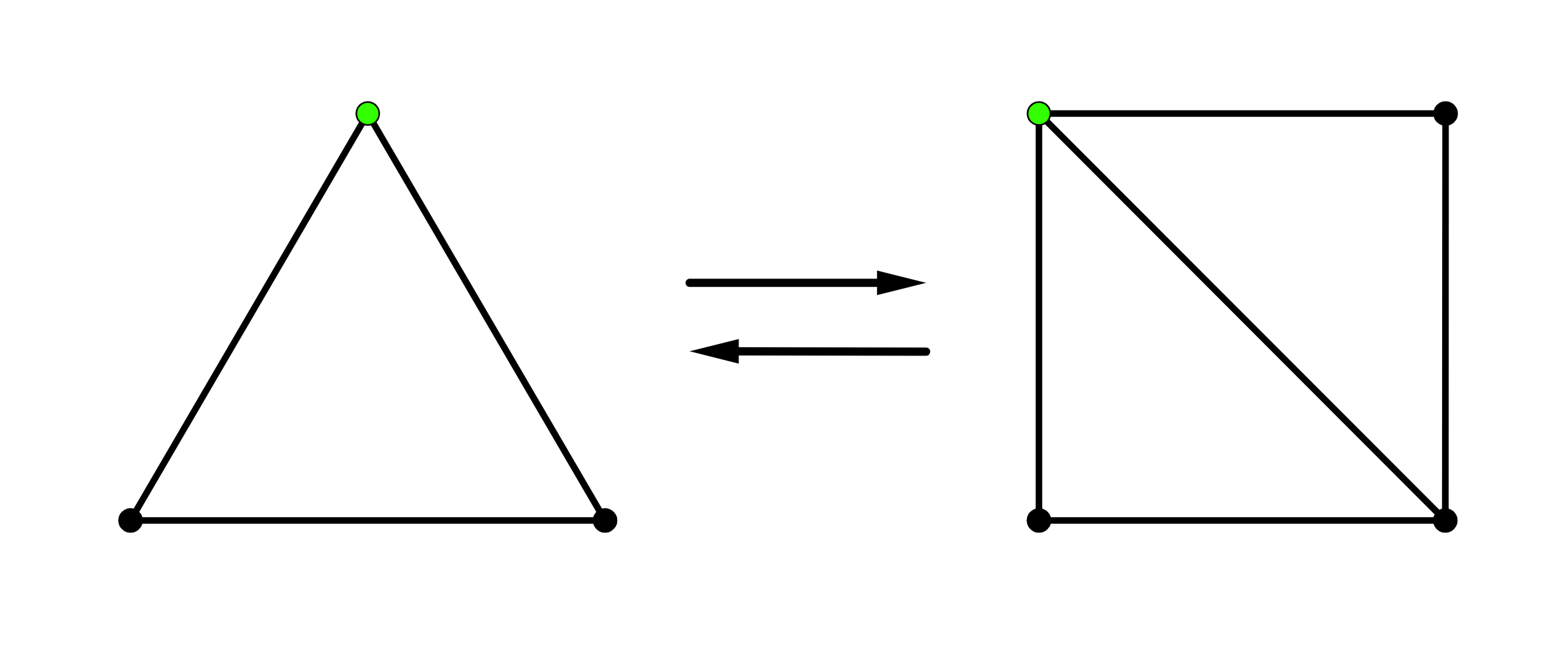}}

	\put(1.02,1.7) {$x$}
	\put(1.82,0.2){$y$}
	\put(0.24,0.2){$z$}

	\put(2.85,1.68){$a$}
	\put(4.32,1.68){$b$}
	\put(4.32,0.2){$c$}
	\put(2.85,0.2){$d$}
	
	\put(2.27, 1.2){$f$}
	\put(2.27, 0.7){$g$}
	\end{picture}
	\caption{Graph Homomorphisms between $C_{3}$ and $R$}
	\label{Graph_Homomorphism_Fig}
\end{figure}
\end{ex}

\begin{df}\cite[Definition 5.1(4)]{BarceloFoundations}\label{Based_Graph_Homomorphism_Def}
A \textit{based graph homomorphism} $f: (X, x_0) \to (Y, y_0)$ is a graph homomorphism $f: X \to Y$ such that $f(x_0) = y_0$.
\end{df}

In Example \ref{Graph_Homomorphism_Ex}, the graph homomorphisms $f$ and $g$ are both based graph homomorphisms. Based simple graphs with loops at every vertex and based graph homomorphisms form a category called $\mathbf{GPH}_{\ast}^{\circ}$.

\begin{df}\cite[p. 74]{HellGraphs}\label{Cartesian_Product_Def}
The \textit{Cartesian product} of the graphs $X$ and $Y$, denoted $X \square Y$, is the graph with vertex set $V(X) \times V(Y)$. There is an edge between the vertices $(x_1, y_1)$ and $(x_2, y_2)$ if either $x_1 = x_2$ and $y_1 y_2 \in E(Y)$ or $y_1 = y_2$ and $x_1 x_2 \in E(X)$. 
\end{df}

The distinguished vertex of the Cartesian product of two based graphs $(X, x_0)$ and $(Y, y_0)$ is the vertex $(x_0, y_0)$. Since there is a loop at every vertex of $X$ and $Y$, there is also a loop at every vertex of $X \square Y$.

\begin{ex}\label{Cartesian_Product_Ex}
The Cartesian product of the graphs $C_3$ and $R$ is illustrated in Figure \ref{Cartesian_Product_Fig}. 
	\begin{figure}[ht]
	\begin{picture}(3.85,2.5)
	\put(0, 0){\includegraphics[scale=0.65]{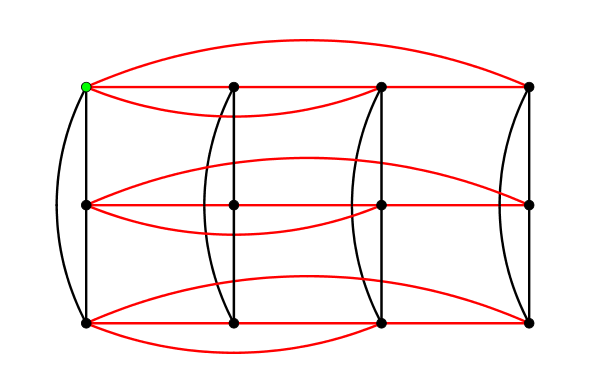}}
		
	\put(0.2, 2.02){$x$}
	\put(0.2, 1.22){$y$}
	\put(0.2, 0.42){$z$}

	\put(0.55,0.08){$a$}
	\put(1.54,0.08){$b$}
	\put(2.53,0.08){$c$}
	\put(3.52,0.08){$d$}
	\end{picture}
	\caption{The Cartesian product of $C_3$ and $R$}
	\label{Cartesian_Product_Fig}
	\end{figure}
The edges of copies of $R$ indexed by vertices of $C_3$ are shown in red. The distinguished vertex of $C_3 \square R$ is $(x, a)$ and shown in green.
\end{ex}

We now proceed to an introduction to A-homotopy theory.

\subsection{A-Homotopy Theory Definitions} 
In algebraic topology, two continuous functions $f, g: A \to B$ are \textit{homotopic} if we can deform $f$ into $g$ over time from 0 to 1 using a continuous function called a \textit{homotopy} \cite[p. 3]{HatcherAlgebraic}. Since the homotopy is a continuous map, it eliminates the need for length and we can do everything over the unit interval. In A-homotopy theory, we deform one graph homomorphism into another graph homomorphism over a path $I_n$ for some $n \in \N$. This length is measured discretely and gives us a combinatorial way to keeps track of the vertices and edges of the graph.

Two spaces $X$ and $Y$ are \textit{homotopy equivalent} if there are continuous maps $f: X \to Y$ and $g: Y \to X$ such that the compositions are homotopic to the identity maps. In particular, we want to know when a space is contractible, or homotopy equivalent to a single point. As stated in the introduction, the \textit{classical fundamental group of a space} is the set of continuous maps from the unit interval into the space under an equivalence. These continuous maps send both endpoints of the interval to the base point of the space. 
The authors of \cite[ Proposition 5.12]{BarceloFoundations}  construct the definition of the homotopy groups of a graph using a cubical set. In this paper, we conceal this cubical set construction to simplify notation. We also only focus on the fundamental group.

\begin{df}\cite[Definition 5.2(1)]{BarceloFoundations}\label{Graph_Homotopy_Def}
Let $f, \; g: (X, x_0) \to (Y, y_0)$ be graph homomorphisms. If there exists a positive integer $n \in \N$ and a graph homomorphism $H: X \square I_{n} \to Y$ such that 
\begin{itemize}
\item $H(v, 0) = f(v)$ for all $v \in V(X)$,
\item $H(v, n) = g(v)$ for all $v \in V(X)$, and
\item $H(x_0, i) = y_0$ for all $0 \leq i \leq n$,
\end{itemize}
then $H$ is a \textit{graph homotopy} from $f$ to $g$. We say that $f$ and $g$ are \textit{A-homotopic}, denoted $f \simeq_{A} g$.
\end{df}

\begin{ex}\label{Graph_Homotopy_Ex}
Recall the graphs $C_3$ and $R$ from Example \ref{Graph_Homomorphism_Ex}. Let $f, g: (C_3, x) \to (R, a)$ be the graph homomorphisms defined by
	\[ \begin{array}{ccc}
	f(x) = a, \;\;\; &			 	& g(x) = a, \\
	f(y) = d, \;\;\; & \textrm{and} \;\;\; & g(y) = b, \\
	f(z) = c, \;\;\; &			 	& g(z) = c.
	\end{array} \]
Figure \ref{f_and_g_Fig} depicts the graph homomorphisms $f$ and $g$. The image under $f$ of each vertex in $C_3$ is shown in red, while the image under $g$ of each vertex in $C_3$ is shown in blue. For example, $f(y) = d$ and $g(y) = b$.	     
    \begin{figure}[ht]
	\begin{picture}(4.6, 1.8)
	\put(0, 0){\includegraphics[scale=0.321]{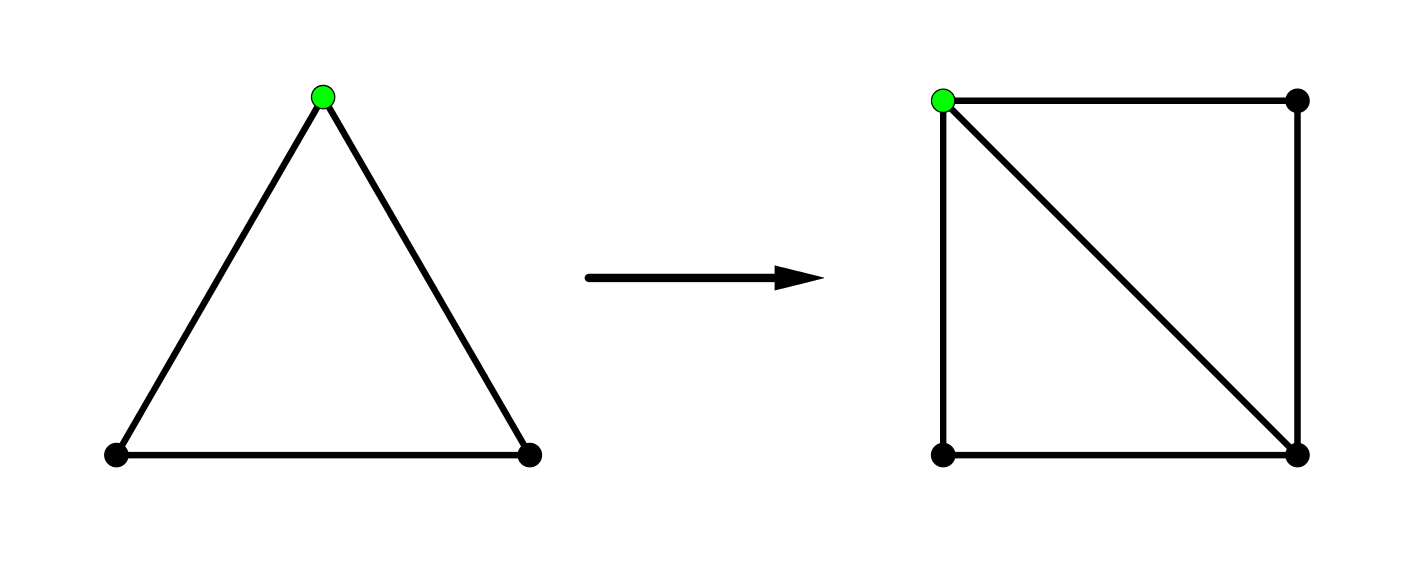}}

	\put(1.03,1.74){$x$}
	\put(1.87,0.26){$y$}
	\put(0.21,0.26){$z$}
	
	\put(2.26,1.15){$f$}
	\put(2.26,0.8){$g$}	
	
	\put(2.99,1.7){$a$}
	\put(4.42,1.7){$b$}
	\put(4.42,0.21){$c$}
	\put(2.99,0.21){$d$}
	
	\put(0.89,1.59){\color{red}{\small $a$}}
	\put(1.68,0.24){\color{red}{\small $d$}}
	\put(0.24,0.45){\color{red}{\small $c$}}
	
	\put(1.18,1.59){\color{blue}{\small $a$}}
	\put(1.85,0.45){\color{blue}{\small $b$}}
	\put(0.41,0.24){\color{blue}{\small $c$}}
	\end{picture}	
	\caption{Graph homomorphisms from $S$ to $T$}
	\label{f_and_g_Fig}
	\end{figure}
	
\noindent Define a map $H: C_3 \square I_{2} \to R$ by 
	\[ \begin{array}{ccc}
	H(x, 0) = a, \;\;\; & H(x, 1) = a, \;\;\; & H(x, 2) = a, \\
	H(y, 0) = d, \;\;\; & H(y, 1) = a, \;\;\; & H(y, 2) = b, \\
	H(z, 0) = c, \;\;\; & H(z, 1) = c, \;\;\; & H(z, 2) = c. 
	\end{array} \]
Figure \ref{Graph_Homotopy_Fig} depicts this map $H$ with the image of each vertex shown in red.
	\begin{figure}[ht]
	\begin{picture}(5, 2.25)
	\put(0, 0){\includegraphics[scale=0.29]{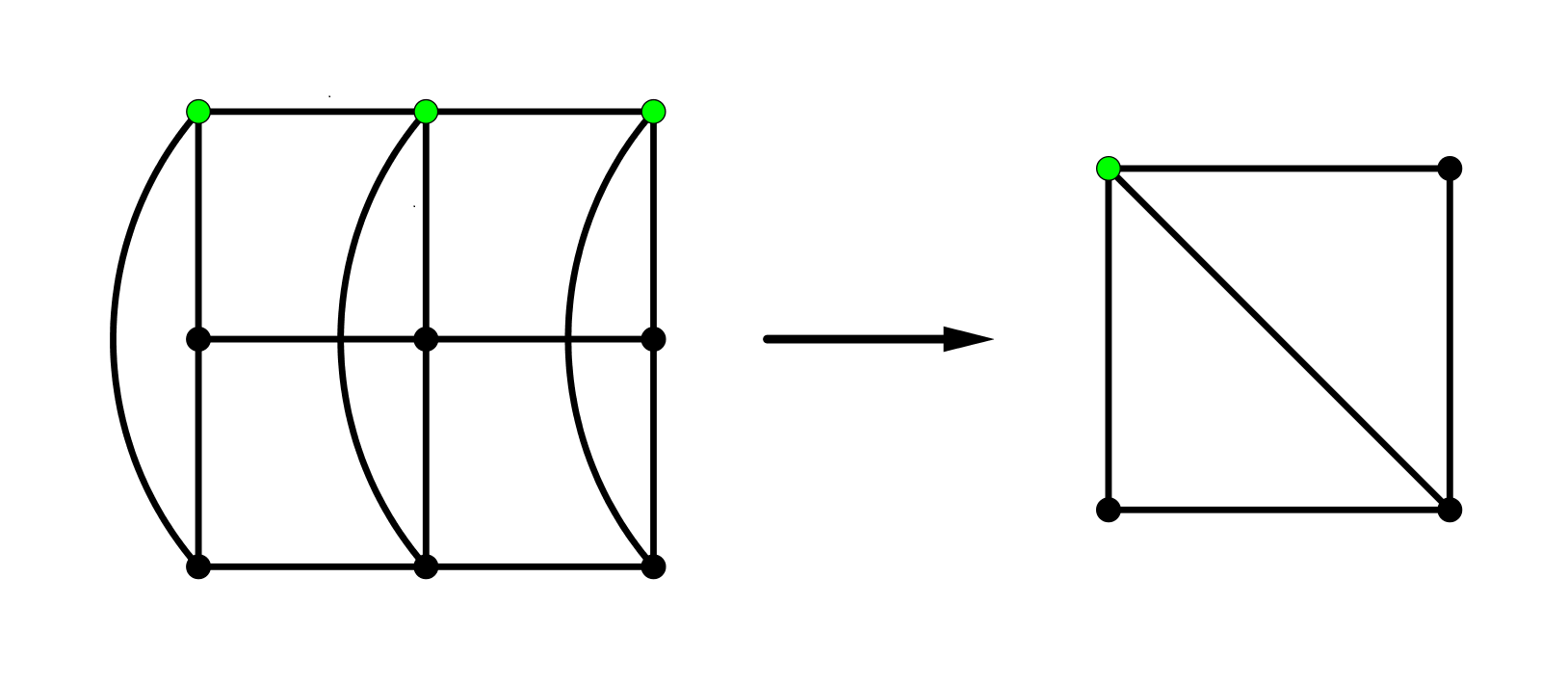}}

	\put(0.2,1.78){$x$}
	\put(0.2,1.05){$y$}
	\put(0.2,0.3){$z$}

	\put(0.57,1.91){$0$}
	\put(1.28,1.91){$1$}
	\put(2,1.91){$2$}
	
	\put(2.65,1.22){$H$}
	
	\put(3.3,1.7){$a$}
	\put(4.65,1.7){$b$}
	\put(4.65,0.4){$c$}
	\put(3.3,0.4){$d$}
	
	\put(0.68,1.67){\color{red}{\small $a$}}
	\put(0.68,0.96){\color{red}{\small $d$}}
	\put(0.68,0.23){\color{red}{\small $c$}}
	
	\put(1.4,1.67){\color{red}{\small $a$}}
	\put(1.4,0.96){\color{red}{\small $a$}}
	\put(1.4,0.23){\color{red}{\small $c$}}
	
	\put(2.12,1.67){\color{red}{\small $a$}}
	\put(2.12,0.96){\color{red}{\small $b$}}
	\put(2.12,0.23){\color{red}{\small $c$}}
	\end{picture}
	\caption{Graph homotopy $H$ from $f$ to $g$} \label{Graph_Homotopy_Fig}
	\end{figure}
Then $H$ is a graph homomorphism with $H(v, 0) = f(v)$ and $H(v, 2) = g(v)$ for all $v \in V(C_3)$, and $H(x, i) = a$ for all $0 \leq i \leq 2$. Thus $H$ is a graph homotopy from $f$ to $g$. However, $H$ is not unique. It is only one of several possible graph homotopies. 
\end{ex}

\begin{df}\cite[Definition 5.2(2)]{BarceloFoundations}\label{Homotopy_Equivalence_Def}
The graph homomorphism $f: X \to Y$ is an \textit{A-homotopy equivalence} if there exists a graph homomorphism $g: Y \to X$ such that $f \circ g \simeq_{A} \mathbf{1}_{Y}$ and $g \circ f \simeq_{A} \mathbf{1}_{X}$. In this case, the graphs $X$ and $Y$ are \textit{A-homotopy equivalent}. 
\end{df}	

\begin{df}\cite{BarceloPerspectives}\label{Concatractible_Def}
A graph $X$ is \textit{A-contractible} if $X$ is A-homotopy equivalent to the graph with a single vertex, called $\ast$, and a single loop edge. For convenience, we will abuse the notation slightly and refer to this graph as $\ast$.
\end{df}

\begin{ex}\cite[p.46]{BarceloPerspectives}\label{C4_Contractible_Ex}
The cycle $C_4$ is A-contractible. 

Let $C_4$ be a $4$-cycle with vertices labelled $a, b, c, d$. There is only one choice for the graph homomorphisms $f$ and $g$. Namely, $f: C_4 \to \ast$ is defined by $f(a) = f(b) = f(c) = f(d) = \ast$ and $g: \ast \to C_4$ is defined by $g(\ast) = a$.

	\begin{figure}[ht]
	\begin{picture}(3.2, 1.9)
	\put(0, 0){\includegraphics[scale=0.5]{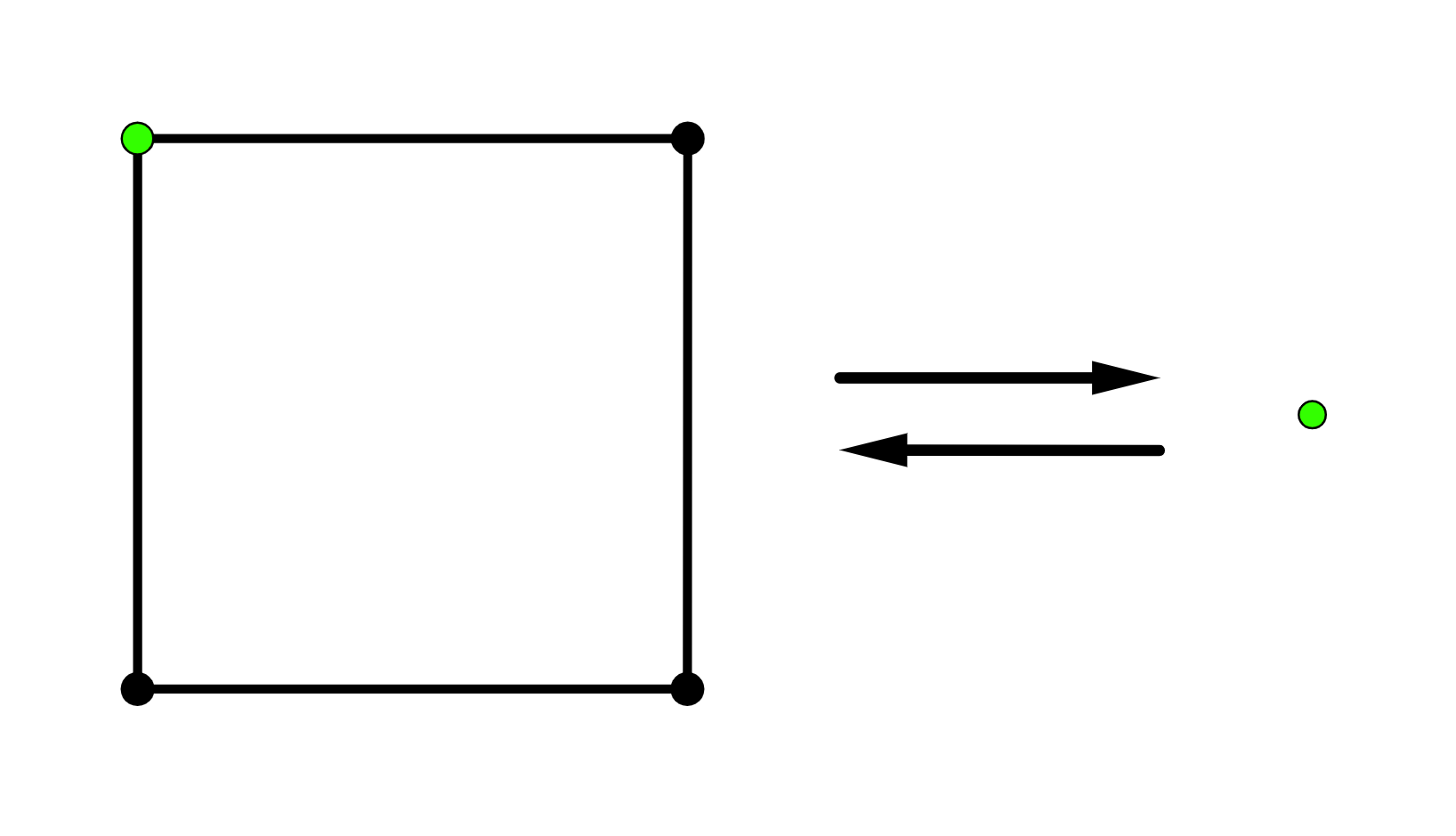}}
		
	\put(0.12,1.6){$a$}
	\put(1.55,1.6){$b$}
	\put(1.55,0.17){$c$}
	\put(0.12,0.17){$d$}
	
	\put(2.95,0.85){{\Large $\ast$}}
	
	\put(2.1,1.1){$f$}
	\put(2.1,0.7){$g$}
	\end{picture}
	\caption{Graph homomorphisms $f$ and $g$}
	\label{C4_Contractible_Fig}
	\end{figure}

\noindent Then $f \circ g$ is defined by $(f \circ g)(\ast) \; = \; f(g(\ast)) \; = \; f(a) \; = \; \ast$, and thus $f \circ g \; = \; \mathbf{1}_{\ast}$.
Also, $g \circ f$ is equal to $c_{a}: C_4 \to C_4$, the constant graph homomorphism mapping every vertex to $a$. To show that $c_{a} \simeq_{A} \mathbf{1}_{C_4}$, define $H: C_4 \square I_{2} \to C_4$ by 
\[ \begin{array}{ccc}
H(a, 0) = a, \;\;\; & H(a, 1) = a, \;\;\; & H(a, 2) = a, \\
H(b, 0) = a, \;\;\; & H(b, 1) = a, \;\;\; & H(b, 2) = b, \\
H(c, 0) = a, \;\;\; & H(c, 1) = d, \;\;\; & H(c, 2) = c, \\
H(d, 0) = a, \;\;\; & H(d, 1) = d, \;\;\; & H(d, 2) = d.
\end{array} \]

	\begin{figure}[ht]
	\begin{picture}(5, 2.5)
	\put(0, 0){\includegraphics[scale=0.48]{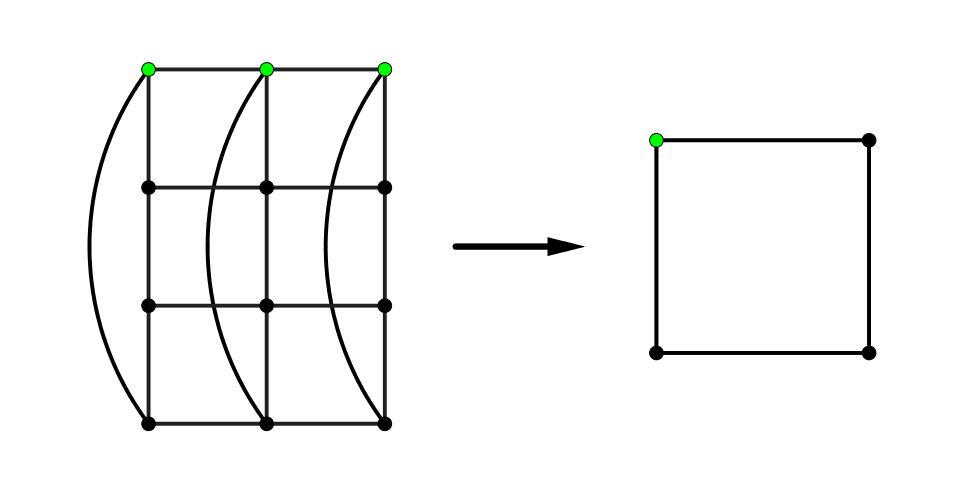}}
		
	\put(0.3,2.1){$a$}
	\put(0.3,1.48){$b$}
	\put(0.3,0.9){$c$}
	\put(0.3,0.3){$d$}
	
	\put(0.68,2.25){$0$}
	\put(1.29,2.25){$1$}
	\put(1.88,2.25){$2$} 
	
	\put(2.45,1.35){$H$}
	
	\put(3.1,1.8){$a$}
	\put(4.45,1.8){$b$}
	\put(4.45,0.58){$c$}
	\put(3.1,0.58){$d$}
	
	\put(0.8,2){\color{red}{\small $a$}}
	\put(0.8,1.4){\color{red}{\small $a$}}
	\put(0.8,0.8){\color{red}{\small $a$}}
	\put(0.8,0.2){\color{red}{\small $a$}}
	
	\put(1.4,2){\color{red}{\small $a$}}
	\put(1.4,1.4){\color{red}{\small $a$}}
	\put(1.4,0.8){\color{red}{\small $d$}}
	\put(1.4,0.2){\color{red}{\small $d$}}
	
	\put(2.02,2){\color{red}{\small $a$}}
	\put(2.02,1.4){\color{red}{\small $b$}}
	\put(2.02,0.8){\color{red}{\small $c$}}
	\put(2.02,0.2){\color{red}{\small $d$}}
	\end{picture}
	\caption{Graph homotopy from $c_{a}$ to $\mathbf{1}_{C_{4}}$}
	\label{Homotopy_C4_Contractible_Fig}
	\end{figure}
	
The image under $H$ of each vertex in $C_4 \square I_{2}$ is shown in red in Figure \ref{Homotopy_C4_Contractible_Fig}. It is routine to verify that $H$ is a graph homomorphism from $c_{a}$ to $ \mathbf{1}_{C_4}$, and $g \circ f \simeq_{A} \mathbf{1}_{C_4}$. Thus the graph $C_4$ is A-contractible.
\end{ex}

In order to construct the fundamental group of a graph, we need to be able to compare graph homomorphisms from paths of any length into the graph. This requires us to work with paths of infinite length.

\begin{nota}\label{n_Fold_Cartesian_Product_Nota}
Let $I_\infty^n$ denote the $n$-fold Cartesian product of $I_\infty$. We will only use non-based graph homomorphisms $f : I_\infty^n
\to X$ with $n = 1, 2$. This will give us paths in a graph $X$ and the graph homotopies between the paths.
\end{nota}

\begin{df}
\cite[Defintion 3.1]{BabsonHomotopy}\label{Stabilizes_Def} Let $f: I_\infty^n \to X$ be a graph homomorphism and $a_1, \ldots , a_n \in \Z$. We say that $f$ \textit{stabilizes in direction $+i$} with $1 \leq i \leq n$, if there is an integer $p^i_f$ such that for all $m \geq p^i_f$, 
$$f(a_1, \ldots , a_{i-1}, m, a_{i+1}, \ldots , a_n) = f(a_1, \ldots , a_{i-1}, p^i_f, a_{i+1}, \ldots , a_n).$$
We say that $f$ \text{stabilizes in the direction $-i$} with $1 \leq i \leq n$ if there is an integer $n^i_f$ such that for all $m \leq n^i_f$,
$$f(a_1, \ldots , a_{i-1}, m, a_{i+1}, \ldots, a_n) = f(a_1, \ldots, a_{i-1}, n^i_f, a_{i+1}, \ldots, a_n).$$
\end{df}

We always take $p^i_f$ to be the least integer and $n^i_f$ to be the greatest integer such that the previous statements are true. If $f$ is constant on the $i^{th}$-axis, then we take $p^i_f = n^i_f = 0$. 

The integers $p^i_f$ and $n^i_f$ give us the points at which the graph homomorphism $f$ stabilizes on the $i^{th}$-axis in the positive and negative directions, respectively. When a graph homomorphism $f : I^n_{\infty} \to X$ stabilizes in every direction, the region of $I^n_\infty$ induced by the vertex set $\prod_{i \in [n]} [n_f^i, p_f^i]$ is called the \textit{active region} of $f$. For each path $f : I_\infty \to X$, we say that $f$ \textit{starts} at $f(n_f^1)$ and $f$ \textit{ends} at $f(p_f^1)$ when these integers exist.

\begin{df}\cite[Defintion 3.1]{BabsonHomotopy} \label{Stable_Graph_Homomorphism_Def} If a graph homomorphism $f : I^n_\infty \to X$ stabilizes in every direction $-i$ and $+i$ for $1 \leq i \leq n$, then we say that $f$ is a \textit{stable graph homomorphism}. Let $S_n(X)$ be the set of stable graph homomorphisms from the infinite $n$-cube $I^n_\infty$ to the graph $X$. 
\end{df}

The set $S_0(X)$ consists of the graph homomorphisms from the graph $\ast$, with a single vertex $\ast$, to the graph $X$. Each graph homomorphism in $S_0(X)$ picks out a vertex of $X$ so that $S_0(X) \cong V(X)$. The graph homomorphisms in $S_1(X)$ give every possible finite walk in $X$.

\begin{ex} \label{Stable_Graph_Homomorphism_Ex}
Figure \ref{Stable_Graph_Homomorphism_Fig} depicts a graph homomorphism $f: I_{\infty} \to C_3$ with the image of each vertex under $f$ shown in red. 
	\begin{figure}[ht]
	\begin{picture}(3.85, 3.8)
	\put(0, 0){\includegraphics[scale=0.365]{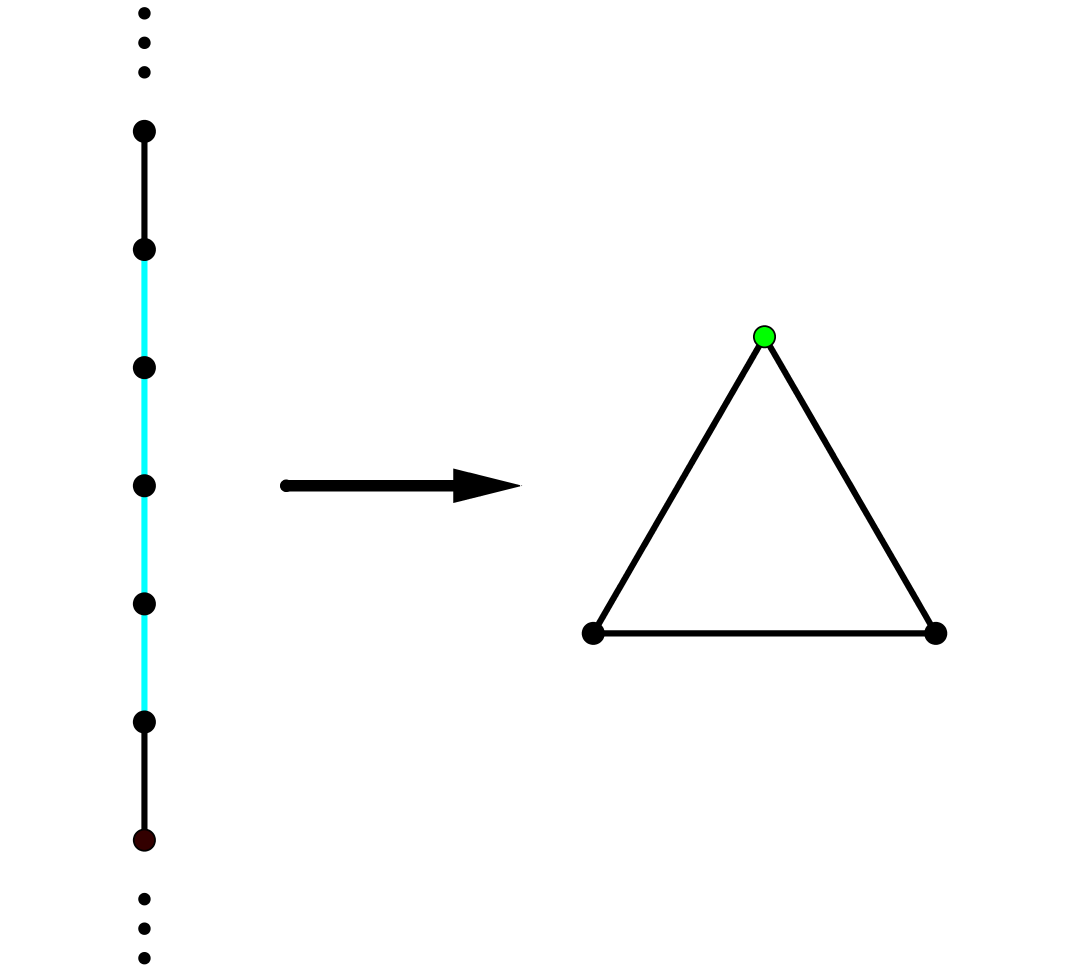}}
		
	\put(0.3,3.15){$4$}
	\put(0.3,2.7){$3$}
	\put(0.3,2.24){$2$}
	\put(0.3,1.79){$1$}
	\put(0.3,1.35){$0$}
	\put(0.23,0.9){$-1$}
	\put(0.23,0.46){$-2$}
	
	\put(0.7,3.15){\color{red}\small{$y$}}
	\put(0.7,2.7){\color{red}\small{$y$}}
	\put(0.7,2.25){\color{red}\small{$x$}}
	\put(0.7,1.79){\color{red}\small{$z$}}
	\put(0.7,1.34){\color{red}\small{$y$}}
	\put(0.7,0.9){\color{red}\small{$x$}}
	\put(0.7,0.46){\color{red}\small{$x$}}
	
	\put(1.4,1.95){$f$}
	
	\put(2.88,2.55){$x$}
	\put(3.6,1.07){$y$}
	\put(2.15,1.07){$z$}
	\end{picture}
	\caption{A stable graph homomorphism $f$ from $I_{\infty}$ to $C_3$}
	\label{Stable_Graph_Homomorphism_Fig}
	\end{figure}
Let $f(i) = x$ for all $i \leq -1$ and  $f(i) = y$ for all $i \geq 3$. Then it follows that the integer $n_f^1 = -1$, that is, $f$ stabilizes on the $1^{st}$-axis in the negative direction at $-1$. Similarly, we have that the integer $p_f^1 = 3$, or $f$ stabilizes on the  $1^{st}$-axis in the positive direction at $3$. The active region of $f$ is shown in blue.
\end{ex}

We can now compare two stable graph homomorphisms using a path homotopy.
 
\begin{df}\cite[Definition 3.2]{BabsonHomotopy}\label{Path_Graph_Homotopy_Def} Let $f, g \in S_1(X)$. The graph homomorphisms $f$ and $g$ are \textit{path homotopic}, denoted $f \sim g$, if there exists a graph homomorphism $H \in S_2(X)$ such that:
\begin{itemize}
    \item[(a)] $f(p_f^1) = g(p_g^1)$ \;\;\; and \;\;\; $f(n_f^1) = g(n_g^1)$,
    \item[(b)] $H(p_H^1, j) = f(p_f^1) = g(p_g^1)$ \;\;\; and \;\;\; \\
    $H(n_H^1, j) = f(n_f^1) = g(n_g^1)$ \;\;\; for all $j \in \Z$,
    \item[(c)] $H(j, n_H^2) = f(j)$ \;\;\; and \;\;\; $H(j, p_H^2) = g(j)$ \;\;\; for all $j \in \Z$.
\end{itemize}
The graph homomorphism $H : I_\infty^2 \to X$ is called a \textit{path graph homotopy from $f$ to $g$}, or \textit{path homotopy} for short.
\end{df}

By part (a), in order for the graph homomorphisms $f, g \in S_1(X)$ to be path homotopic, they must stabilize to the same vertex in the direction -1 and to the same vertex in the direction +1. By part (b), the graph homomorphism $H$ must stabilize in the directions -1 and +1 to the same vertices that $f$ and $g$ stabilize to in the directions -1 and +1 repeated along the $1^{st}$-axis. By part (c), the graph homomorphism $H$ must stabilize to $f$ in the direction -2 and stabilize to $g$ in the direction +2.

\begin{ex}\label{Path_Graph_Homotopy_Ex}
Recall the graph $R$ from Ex \ref{Graph_Homomorphism_Ex}. Let the graph homomorphisms $f, g \in S_1(R)$ be defined by
$$f(i) = \begin{cases}
a & \text{for} \;\;\; i \geq 2, \\
c & \text{for} \;\;\; i = 1, \\
d & \text{for} \;\;\; i = 0, \\
a & \text{for} \;\;\; i \leq -1,
\end{cases} 
\;\;\;\;\; \text{and} \;\;\;\;\; 
g(i) = \begin{cases}
a & \text{for} \;\;\; i \geq 2, \\
c & \text{for} \;\;\; i = 1, \\
b & \text{for} \;\;\; i = 0, \\
a & \text{for} \;\;\; i \leq -1.
\end{cases} $$
Consider the graph homomorphism $H \in S_2(R)$ defined in Figure \ref{Path_Graph_Homotopy_Fig}.

	\begin{figure}[ht]
	\begin{picture}(4.5, 3)
	\put(0, 0){\includegraphics[scale=0.17]{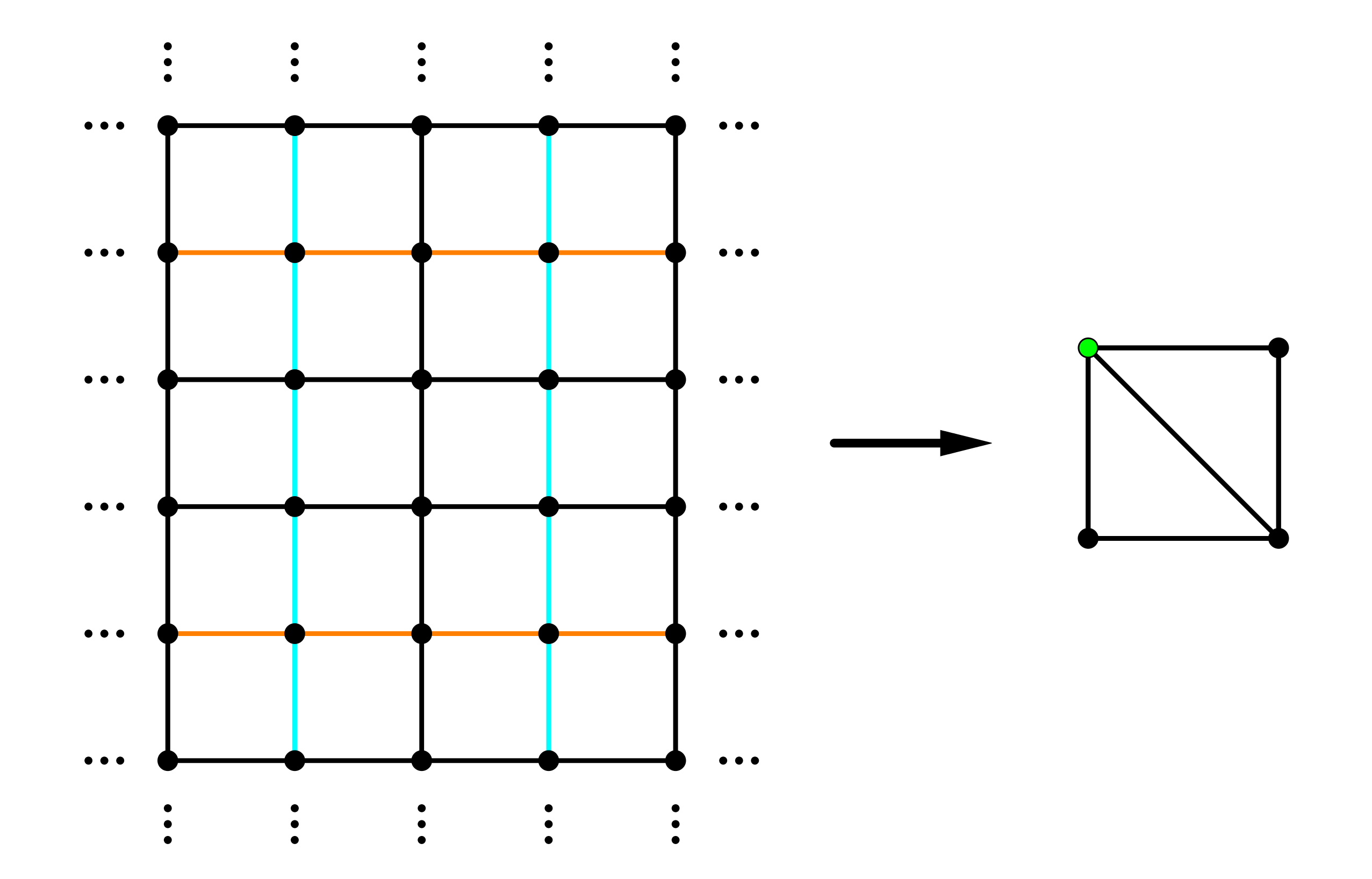}}

	\put(0.1,2.43){$3$}
	\put(0.1,2.02){$2$}
	\put(0.1,1.6){$1$}
	\put(0.1,1.17){$0$}
	\put(0,0.77){$-1$}
	\put(0,0.35){$-2$}

	\put(0.38,2.85){$-1$}
	\put(0.92,2.85){$0$}
	\put(1.33,2.85){$1$}
	\put(1.75,2.85){$2$}
	\put(2.2,2.85){$3$}
	
	\put(2.9, 1.5){$H$}
	
	\put(3.45,1.85){$a$}
	\put(4.28,1.85){$b$}
	\put(4.28,0.95){$c$}
	\put(3.45,0.95){$d$}
	
	\put(0.59,2.38){\color{red}\small{$a$}}
	\put(0.59,1.96){\color{red}\small{$a$}}
	\put(0.59,1.54){\color{red}\small{$c$}}
	\put(0.59,1.1){\color{red}\small{$d$}}
	\put(0.59,0.7){\color{red}\small{$a$}}
	\put(0.59,0.28){\color{red}\small{$a$}}
	
	\put(1.01,2.38){\color{red}\small{$a$}}
	\put(1.01,1.96){\color{red}\small{$a$}}
	\put(1.01,1.54){\color{red}\small{$c$}}
	\put(1.01,1.1){\color{red}\small{$d$}}
	\put(1.01,0.7){\color{red}\small{$a$}}
	\put(1.01,0.28){\color{red}\small{$a$}}
	
	\put(1.44,2.38){\color{red}\small{$a$}}
	\put(1.44,1.96){\color{red}\small{$a$}}
	\put(1.44,1.54){\color{red}\small{$c$}}
	\put(1.44,1.1){\color{red}\small{$c$}}
	\put(1.44,0.7){\color{red}\small{$a$}}
	\put(1.44,0.28){\color{red}\small{$a$}}
	
	\put(1.85,2.38){\color{red}\small{$a$}}
	\put(1.85,1.96){\color{red}\small{$a$}}
	\put(1.85,1.54){\color{red}\small{$c$}}
	\put(1.85,1.1){\color{red}\small{$b$}}
	\put(1.85,0.7){\color{red}\small{$a$}}
	\put(1.85,0.28){\color{red}\small{$a$}}
	
	\put(2.28,2.38){\color{red}\small{$a$}}
	\put(2.28,1.96){\color{red}\small{$a$}}
	\put(2.28,1.54){\color{red}\small{$c$}}
	\put(2.28,1.1){\color{red}\small{$b$}}
	\put(2.28,0.7){\color{red}\small{$a$}}
	\put(2.28,0.28){\color{red}\small{$a$}}
	\end{picture}
	\caption{A stable graph homomorphism $H$ from $I_\infty^2$ to $R$}
	\label{Path_Graph_Homotopy_Fig}
	\end{figure}

For the purpose of illustration, here is how to verify that $H$ is a path homotopy from $f$ to $g$ as in Definition \ref{Path_Graph_Homotopy_Def}.
\begin{itemize}
    \item[(a)] By definition of $f$ and $g$, we have $p_f^1 = p_g^1 = 2$ and $n_f^1 = n_g^1 = -1$. In particular, $f(2) = a = g(2)$ and $f(-1) = a = g(-1)$. Thus $f$ and $g$ both stabilize to the vertex $a$ in the directions $+1$ and $-1$.
    \item[(b)] By definition of $H$, we have $p_H^1 = 2$ and $n_H^1 = -1$. The areas of $I_\infty^2$ where $H$ stabilizes on the $1^{st}$-axis are shown in orange in Figure \ref{Path_Graph_Homotopy_Fig}. For all $j \in \Z$, $H(2, j) = a = f(2) = g(2)$, and $H$ stabilizes in the $+1$ direction to the vertex $a$ repeated along the $1^{st}$-axis. Similarly, $H(-1, j) = a = f(-1) = g(-1)$ for all $j \in \Z$, and $H$ stabilizes in the $-1$ direction to the vertex $a$ repeated along the $1^{st}$-axis.  
    \item[(c)] By definition of $H$, we have $n_H^2 = 0$ and $p_H^2 = 2$. The areas of $I_\infty^2$ where $H$ stabilizes on the $2^{st}$-axis are shown in blue in Figure \ref{Path_Graph_Homotopy_Fig}. For all $j \in \Z$, $H(j, 0) = f(j)$ and $H(j, 2) = g(j)$, and $H$ stabilizes to $f$ in the direction $-2$ and stabilizes to $g$ in the direction $+2$.
\end{itemize}
Thus $H$ is a graph homotopy from $f$ to $g$, and $f \sim g$.
\end{ex}

This path homotopy relation $\sim$ gives an equivalence relation on $S_1(X)$ \cite[Proposition 3.3]{BabsonHomotopy}.

\begin{df} \cite[Definition 3.4]{BabsonHomotopy}\label{Based_Stable_Set_Def}
Let $x_0 \in X$ be a distinguished vertex of the graph $X$. The set $B_1(X, x_0) \subset S_1(X)$ is the subset of all graph homomorphisms from $I_\infty$ to $X$ that stabilize to $x_0$ in the directions $-1$ and $+1$.
\end{df}

\begin{df} \cite[Proposition 3.5]{BabsonHomotopy}\label{Fundamental_Group_Def} Let $A_1(X, x_0) = B_1(X, x_0)/ \sim$. The set $A_1(X, x_0)$ is the \textit{fundamental group} of the graph $(X, x_0)$.
\end{df}

By Proposition 3.5 of \cite{BabsonHomotopy}, this construction is a group with the operation of concatenation, which we define in Definition \ref{Concatenation_Def}. The identity of this group is the equivalence class of the path $p_{x_0}: I_\infty \to X$ given by $p_{x_0}(i) = x_0$ for all $i \in \Z$. The inverse of the equivalence class of a path $\gamma : I_\infty \to X$ is the equivalence class of the path $\overline{\gamma}: I_\infty \to X$ given by $\overline{\gamma}(i) = \gamma(-i)$ for all $i \in \Z$. We now explore the group operation on $B_1(X, x_0)/ \sim$ and its properties.

\begin{df}\cite[p. 34]{BabsonHomotopy}\label{Concatenation_Def}
Let $f$ and $g$ be graph homomorphisms of $S_1(X)$ with $f(n_f^1) = g(p_g^1)$. The \textit{concatenation of $f$ and $g$}, denoted $f \cdot g$, is defined by
$$(f \cdot g)(i) = \begin{cases}
f(i + n_f^1) & \text{for} \;\;\; i \geq 0, \\
g(i + p_f^1) & \text{for} \;\;\; i \leq 0.
\end{cases}$$
\end{df}

This operation essentially shifts the first graph homomorphism $f$ to stabilize in the direction $-1$ at zero and shifts the second graph homomorphism $g$ to stabilize in the direction $+1$ at zero. For this reason, $f$ must stabilize in the direction $-1$ to the same vertex that $g$ stabilizes to in the direction $+1$. The following result will be utilized throughout the rest of the paper.

\begin{lemma}\label{Concatenation_Prop}
If $f, g \in S_1(X)$ with $f(n_f^1) = g(p_g^1)$, then the concatenation $f \cdot g$ is a graph homomorphism of $S_1(X)$ that stabilizes in the direction $+1$ at $p_{f \cdot g}^1 = p_f^1 - n_f^1$ and in the direction $-1$ at $n_{f \cdot g}^1 = n_g^1 - p_g^1$.
\end{lemma}

\begin{proof}
Let $f, g \in S_1(X)$ with $f(n_f^1) = g(p_g^1)$. The concatenation $f \cdot g$ is well-defined, since $f(n_f^1) = g(p_g^1)$. In order for $f \cdot g$ to be a graph homomorphism, each pair of adjacent vertices in $I_\infty$ must be mapped to adjacent vertices in $X$. By definition of the graph $I_\infty$, there are edges $i i, i(i + 1) \in E(I_\infty)$ for each $i \in \Z$.
\begin{itemize}
    \item If $i \geq 0$, then 
    $$(f \cdot g)(i) = f(i + n_f^1) \;\;\; \text{and} \;\;\;
    (f \cdot g)(i + 1) = f(i + 1 + n_f^1).$$
    \item Otherwise $i < 0$, and it follows that
    $$(f \cdot g)(i) = g(i + p_g^1)
    \;\;\; \text{and} \;\;\; 
    (f \cdot g)(i + 1) = g(i + 1 + p_g^1).$$
\end{itemize}
These pairs of vertices must be adjacent, since $f$ and $g$ is a graph homomorphism. Therefore, the concatenation $f \cdot g$ is a graph homomorphism.

Now we must examine where $f \cdot g$ stabilizes in direction $-1$ and $+1$. For $i \geq 0$, $(f \cdot g)(i) = f(i + n_f^1)$. Since $p_f^1 - n_f^1 \geq 0$, 
$$(f \cdot g)(p_f^1 - n_f^1) = f(p_f^1 - n_f^1 + n_f^1) = f(p_f^1).$$
By Definition \ref{Stabilizes_Def}, $p_f^1$ is the least integer such that $f(m) = f(p_f^1)$ for all $m \geq p_f^1$, so it follows that $p_f^1 - n_f^1$ is the least integer such that $(f \cdot g)(i) = f(p_f^1)$ for all $i \geq p_f^1 - n_f^1$. Therefore, $p_{f \cdot g}^1 = p_f^1 - n_f^1$. Similarly, we have that $n_{f \cdot g}^1 = n_g^1 - p_g^1$. 
\end{proof}

\subsection{Properties of A-Homotopy}
Now we explore some properties of the A-homotopy fundamental group and graph homomorphisms. These properties will be useful in Section 3 to prove the lifting properties of covering graphs.

\begin{df}\label{Induced_Map_Def}
Let $f : (Y, y_0) \to (X, x_0)$ be a graph homomorphism. The \textit{induced map} $f_\ast : A_1(Y, y_0) \to A_1(X, x_0)$ is defined by $f_\ast([\gamma]) = [f \circ \gamma]$, where $[\gamma]$ is an equivalence class of $A_1(Y, y_0)$.
\end{df}

\begin{lemma}\label{Induced_Map_Well_Defined_Lemma}
If $f : (Y, y_0) \to (X, x_0)$ is a graph homomorphism, then the induced map $f_\ast: A_1(Y, y_0) \to A_1(X, x_0)$ is well-defined.
\end{lemma}

\begin{proof}
Let $f : (Y, y_0) \to (X, x_0)$ is a graph homomorphism and $\gamma_1, \gamma_2 \in B_1(Y, y_0)$ such that $\gamma_1 \sim \gamma_2$. Then there exist a path homotopy $H \in S_2(Y)$ from $\gamma_1$ to $\gamma_2$. 
Define $H': I_\infty^2 \to X$ by $H' = f \circ H$. Since $H'$ is a composition of graph homomorphisms, it follows that $H'$ is a graph homomorphism. The conditions of Definition \ref{Path_Graph_Homotopy_Def} are preserved by composition with $f$, so $H'$ is a path homotopy from $f \circ \gamma_1$ to $f \circ \gamma_2$. Thus $f_\ast$ is well-defined. 
\end{proof}

The two next lemmas describe two cases for when stable graph homomorphisms are path homotopic. The Shifting Lemma (\ref{Shifting_Lemma}) states that a path is homotopic to that same path shifted to start at a different vertex.

\begin{lemma}[Shifting Lemma] \label{Shifting_Lemma}
Let $f: I_\infty \to X$ be a stable graph homomorphism and $n \in \Z$. Suppose $f': I_\infty \to X$ be a stable graph homomorphism such that $f'(i) = f(i - n)$, that is, $f$ shifted by $n$. Then $f \sim f'$.
\end{lemma}

\begin{proof}
Let $f, f'$ be as stated in the lemma. Then $f \sim f'$ by the graph homotopy $H': I_\infty^2 \to X$ defined by
$$H'(i, j) = 
\begin{cases}
f(i) & \text{for} \;\;\; j \leq 0, \\
f(i-j) & \text{for} \;\;\; 0 \leq j \leq n, \\
f(i-n) & \text{for} \;\;\; j \geq n,
\end{cases}$$
if $n \geq 0$ and by
$$H'(i, j) = 
\begin{cases}
f(i) & \text{for} \;\;\; j \leq 0, \\
f(i+j) & \text{for} \;\;\; 0 \leq j \leq -n, \\
f(i-n) & \text{for} \;\;\; j \geq -n,
\end{cases}$$
if $n \leq 0$.
\end{proof}

When a path $f : I_\infty \to X$ maps a sequence of consecutive vertices to the same vertex in $X$, we call this section \textit{padding}. The Padding Lemma (\ref{Padding_Lemma}) states that a path with padding is homotopic to that same path with the padding removed.

\begin{lemma}[Padding Lemma]\label{Padding_Lemma}
Let $f: I_\infty \to X$ be a stable graph homomorphism. Suppose $f': I_\infty \to X$ is a stable graph homomorphism such that
$$f'(i) = \begin{cases}
f(i - m) & \text{ for} \;\;\; i \geq b + m, \\
f(b) & \text{for} \;\;\; b - n \leq i \leq b + m, \\
f(i + n) & \text{for} \;\;\; i \leq b - n,
\end{cases}$$
for some $n, m \in \N$ and $b \in \Z$ with $n_f^1 < b < p_f^1$. Then $f \sim f'$.
\end{lemma}

\begin{proof}
Let $f, f'$ be as stated in the lemma. Define $g: I_\infty \to X$ by $g(i) = f'(i-n)$ for all $i \in \Z$.
Then $f \sim g$ by the path homotopy $H: I_\infty^2 \to X$ given by 
$$H(i, j) = \begin{cases}
f(i) & \text{for} \;\;\; j \leq 0, \\
f(i - j) & \text{for} \;\;\; 0 \leq j \leq n + m, \; i \geq b + j, \\
f(b) & \text{for} \;\;\; 0 \leq j \leq n + m, \; b \leq i \leq b + j, \\
f(i) & \text{for} \;\;\; 0 \leq j \leq n + m, \; i \leq b, \\
g(i) & \text{for} \;\;\; j \geq n + m.
\end{cases}$$
Since $g$ is $f'$ shifted by $n$, it follows that $g \sim f'$ by the Shifting Lemma \ref{Shifting_Lemma}. Since the path homotopy relation $\sim$ is an equivalence relation, this implies that $f \sim f'$.
\end{proof}

We use the Shifting Lemma (\ref{Shifting_Lemma}) in the proofs of Lemma \ref{Gamma_n_Inverse_Lemma} and Proposition \ref{Fundamental_Group_Ck_Thm}, and we use the Padding Lemma (\ref{Padding_Lemma}) in the proofs of Lemmas \ref{Induced_Map_Group_Homomorphism_Lemma} and \ref{Gamma_n_Lift_Equivalence_Class_Lemma}. The following lemma is useful in proofing the Lifting Criterion (Theorem \ref{Lifting_Criterion}).

\begin{lemma}\label{Induced_Map_Group_Homomorphism_Lemma}
If $f : (Y, y_0) \to (X, x_0)$ is a graph homomorphism, then the induced map $f_\ast: A_1(Y, y_0) \; \to A_1(X, x_0)$ is a group homomorphism.
\end{lemma}

\begin{proof}
Let  $f : (Y, y_0) \to (X, x_0)$ is a graph homomorphism and $\gamma_1, \gamma_2 \in B_1(Y, y_0)$. Since $B_1(Y, y_0)$ is closed with respect to concatenation, it follows that $\gamma_1 \cdot \gamma_2 \in B_1(Y, y_0)$. We must show that $f_\ast([\gamma_1 \cdot \gamma_2]) = f_\ast([\gamma_1]) \cdot f_\ast([\gamma_2])$. 
The composition $f \circ (\gamma_1 \cdot \gamma_2)$ is defined by
$$f \circ (\gamma_1 \cdot \gamma_2)(i) = \begin{cases}
f(\gamma_1(i + n_{\gamma_1}^1)) & \text{for} \;\;\; i \geq 0, \\
f(\gamma_2(i + p_{\gamma_2}^1)) & \text{for} \;\;\; i \leq 0.
\end{cases}$$
Similarly, the concatenation $(f \circ \gamma_1) \cdot (f \circ \gamma_2)$ is defined by
\begin{eqnarray*}
((f \circ \gamma_1) \cdot (f \circ \gamma_2))(i)
& = & \begin{cases}
f(\gamma_1(i + n_{f \circ \gamma_1}^1)) & \text{for} \;\;\; i \geq 0, \\
f(\gamma_2(i + p_{f \circ \gamma_2}^1)) & \text{for} \;\;\; i \leq 0.
\end{cases}
\end{eqnarray*}
We need to relate $n_{\gamma_1}^1$ to $n_{f \circ \gamma_1}^1$ and $p_{\gamma_2}^1$ to $p_{f \circ \gamma_2}^1$. Since $f: (Y, y_0) \to (X, x_0)$ might map some vertices of $Y$ to the same vertex in $X$, it follows that $n_{f \circ \gamma_1}^1 \geq n_{\gamma_1}^1$
and $p_{f \circ \gamma_2}^1 \leq p_{\gamma_2}^1$. Thus $n_{f \circ \gamma_1}^1 - n_{\gamma_1}^1 \geq 0$ and $p_{\gamma_2}^1 - p_{f \circ \gamma_2}^1 \geq 0$. By adding zero, we have that
\begin{eqnarray*}
((f \circ \gamma_1) \cdot (f \circ \gamma_2))(i)
& = & \begin{cases}
f(\gamma_1(i + n_{f \circ \gamma_1}^1 + n_{\gamma_1}^1 - n_{\gamma_1}^1)) & \text{for} \;\;\; i \geq 0, \\
f(\gamma_2(i + p_{f \circ \gamma_2}^1 + p_{\gamma_2}^1 - p_{\gamma_2}^1)) & \text{for} \;\;\; i \leq 0,
\end{cases} \\
& = & \begin{cases}
f(\gamma_1 \cdot \gamma_2)(i + (n_{f \circ \gamma_1}^1  - n_{\gamma_1}^1)) & \text{for} \;\;\; i \geq 0, \\
f(\gamma_1 \cdot \gamma_2)(i - (p_{\gamma_2}^1 - p_{f \circ \gamma_2}^1)) & \text{for} \;\;\; i \leq 0.
\end{cases} 
\end{eqnarray*}
This shows that aside from some potentially padding in $(f \circ \gamma_1) \cdot (f \circ \gamma_2)$  from the vertex $p_{f \circ \gamma_2}^1 - p_{\gamma_2}^1$ to the vertex $n_{f \circ \gamma_1}^1 - n_{\gamma_1}^1$, it is the same as $f \circ (\gamma_1 \cdot \gamma_2)$. Therefore, $f \circ (\gamma_1 \cdot \gamma_2) \sim (f \circ \gamma_1) \cdot (f \circ \gamma_2)$ by the Padding Lemma \ref{Padding_Lemma}, and it follows that $f_\ast$ is a group homomorphism.
\end{proof}

\begin{lemma}\label{Induced_Map_Composition_Lemma}
Let $g: (Y,y_0) \to (Z, z_0)$, and $h: (Z, z_0) \to (X, x_0)$ be graph homomorphisms such that $h \circ g = f$. Then $h_\ast \circ g_\ast = (h \circ g)_\ast$.
\end{lemma}

This result follows directly from Definition \ref{Induced_Map_Def} and the associativity of composition. This result implies that $A_1: \mathbf{GPH}_{\ast}^{\circ} \to \mathbf{GRP}$ is a functor with $A_1(f) = f_\ast$.

\section{Lifting Properties}
This section contains the main contributions of this paper, the lifting properties of covering graphs. A-homotopy theory is a homotopy theory for graphs; however, much of the surrounding theory has not been developed yet. Here, we develop some of the analogous covering space theory by following the approach in \cite{HatcherAlgebraic} as much as possible. In topology, a \textit{covering space} of a space $A$ is a space $\widetilde{A}$ with a continuous map $\rho: \widetilde{A} \to A$ that preserves the local structure of the space. When considering a graph as a subspace of $\R^3$, these covering spaces fail to recognize the structure of the graph, namely, the vertices and edges. Thus we require covering graphs of a graph $X$, that is, graphs $\tX$ with graph homomorphisms $\rho: \tX \to X$ that preserve the local structures of the graphs. Given a covering space $(\widetilde{A}, \rho)$ of $A$ and a continuous map $f: B \to A$, a \textit{lift} of $f$ is a map $\widetilde{f}: B \to \widetilde{A}$ which factors $f$ through the space $\widetilde{A}$. In topology, there are certain properties that determine when a lift does or does not exist for covering spaces, called lifting properties. We will define a discrete version of lifts and develop the corresponding lifting properties in this section. 

The next three definitions give us a more precise idea of what covering graphs are. First, we introduce a vertex set analogous to the closed and open sets used in topology.

\begin{df}\label{Closed_Neighborhood_Def}
Let $X$ be a graph, and let $x \in V(X)$. The \textit{closed neighborhood of $x$} is the set of vertices
$$N[x] = \{v \in V(X) \;| \; vx \in E(X)\}.$$
\end{df}

Since $X$ has a loop at every vertex, we have that $x \in N[x]$ for all $x \in V(X)$.The following definition is the main component of a covering graph and standard in graph theory. It can be found in \cite{GodsilAlgebraic}.

\begin{df}\cite[p. 115]{GodsilAlgebraic}\label{Local_Isomorphism_Def}
The graph homomorphism $\rho: X \to Y$ is a \textit{local isomorphism} if for each vertex $y \in V(Y)$ and each vertex $x \in \rho^{-1}(y)$, the vertex set map restricted to $N[x]$, $\rho|_{N[x]}: N[x] \to N[y]$, is bijective.
\end{df}

One important property of local isomorphisms is that they are locally invertible. This property is useful in proving the lifting properties in Section 3.  

\begin{df}\label{Neighborhood_Subgaph_Def}
Let $X$ be a graph. For $x \in V(X)$, the \textit{neighborhood subgraph of $x$}, denoted $N_x$, is the subgraph of $X$ with vertex set $V(N_x) = N[x]$ and edge set $E(N_x) = \{ xv \;|\; v \in N[x]\} \cup \{ v v \; | \; v \in N[x] \}$.
\end{df}

If $\rho: X \to Y$ is a local isomorphism, then $\rho$ induces a graph homomorphism from the subgraph $N_x$ to the subgraph $N_{\rho(x)}$ for each $x \in V(X)$. Thus there is a graph homomorphism $\rho|_{N_x}: N_x \to N_{\rho(x)}$, which is bijective on the vertices and edges of the subgraphs. This implies the following lemma.

\begin{lemma}\label{Relative_Local_Isomorphism_Lemma}
Let $\rho: X \to Y$ be a local isomorphism and $x \in V(X)$. Then the graph homomorphism $\rho|_{N_x}$ is invertible, and its inverse $(\rho|_{N_x})^{-1}: N_{\rho(x)} \to N_x$ is a graph homomorphism.
\end{lemma}

For a graph $X$ and a vertex $x \in V(X)$, the neighborhood subgraph $N_x$ should not be confused with the induced subgraph of the vertex set $N[x]$. These subgraphs are not necessarily equal.

\begin{df}\cite[p. 3]{GodsilAlgebraic}\label{Induced_Subgraph_Def}
Let $X$ be a graph and $V' \subseteq V(X)$. The \textit{induced subgraph} $X[V']$ is the graph with vertex set $V'$ and edge set $E' = \{ vw \in E(X) \;|\; v, w \in V'\}$.
\end{df}

For a local isomorphism $\rho: X \to Y$, the graph homomorphism $\rho|_{X[N[x]]}: X[N[x]] \to X[N[\rho(x)]]$, induced by $\rho$, is not always bijective on the edges of the subgraphs. Thus this graph homomorphism $\rho|_{X[N[x]]}$ is not necessarily invertible locally.

\begin{ex}\label{Local_Isomorphism_Ex}
Let $\rho: C_6 \to C_3$ be the local isomorphism with $\rho([i]) = [i \mod 3]$ for $i \in \{0, \ldots , 6\}$ depicted in Figure \ref{Local_Iso}. On the left, the edges of the neighborhood subgraphs $N_{[4]}$ in $C_6$ and $N_{[1]}$ in $C_3$ are shown in blue. The graph homomorphism $\rho|_{N_{[4]}}: N_{[4]} \to N_{[1]}$ is invertible. On the right, the edges of the induced subgraphs $C_6[N[4]]$ and $C_3[N[1]]$ are highlighted in blue. Since there is an edge $[0][2]$ in the subgraph $C_3[N[1]]$ and no edge $[3][5]$ in the subgraph $C_6[N[4]]$, the graph homomorphism $\rho|_{C_6(N[4])}$ is not invertible.

\begin{figure}[ht] 
	\begin{picture}(3.9,3.2)
	\put(0,0){\includegraphics[scale=0.48]{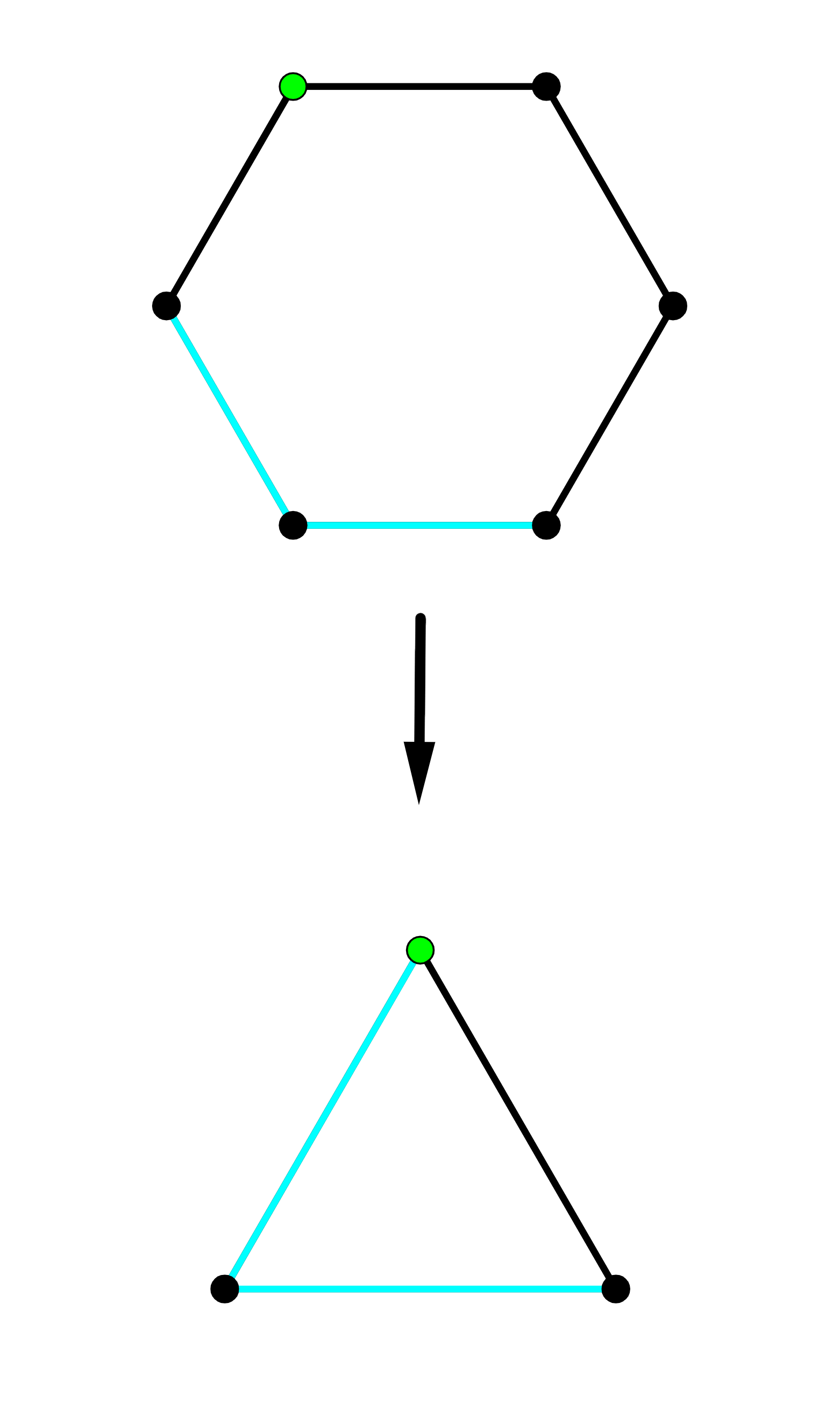}}
	
	\put(0.48,3.02){{\footnotesize $[0]$}}
	\put(1.25,3.02){{\footnotesize $[1]$}}
	\put(1.55,2.42){{\footnotesize $[2]$}}
	\put(1.25,1.82){{\footnotesize $[3]$}}
	\put(0.48,1.82){{\footnotesize $[4]$}}
	\put(0.19,2.42){{\footnotesize $[5]$}}

	\put(0.87,1.1){{\footnotesize $[0]$}}
	\put(0.33,0.2){{\footnotesize $[1]$}}
	\put(1.42,0.2){{\footnotesize $[2]$}}	

	\put(2,0){\includegraphics[scale=0.48]{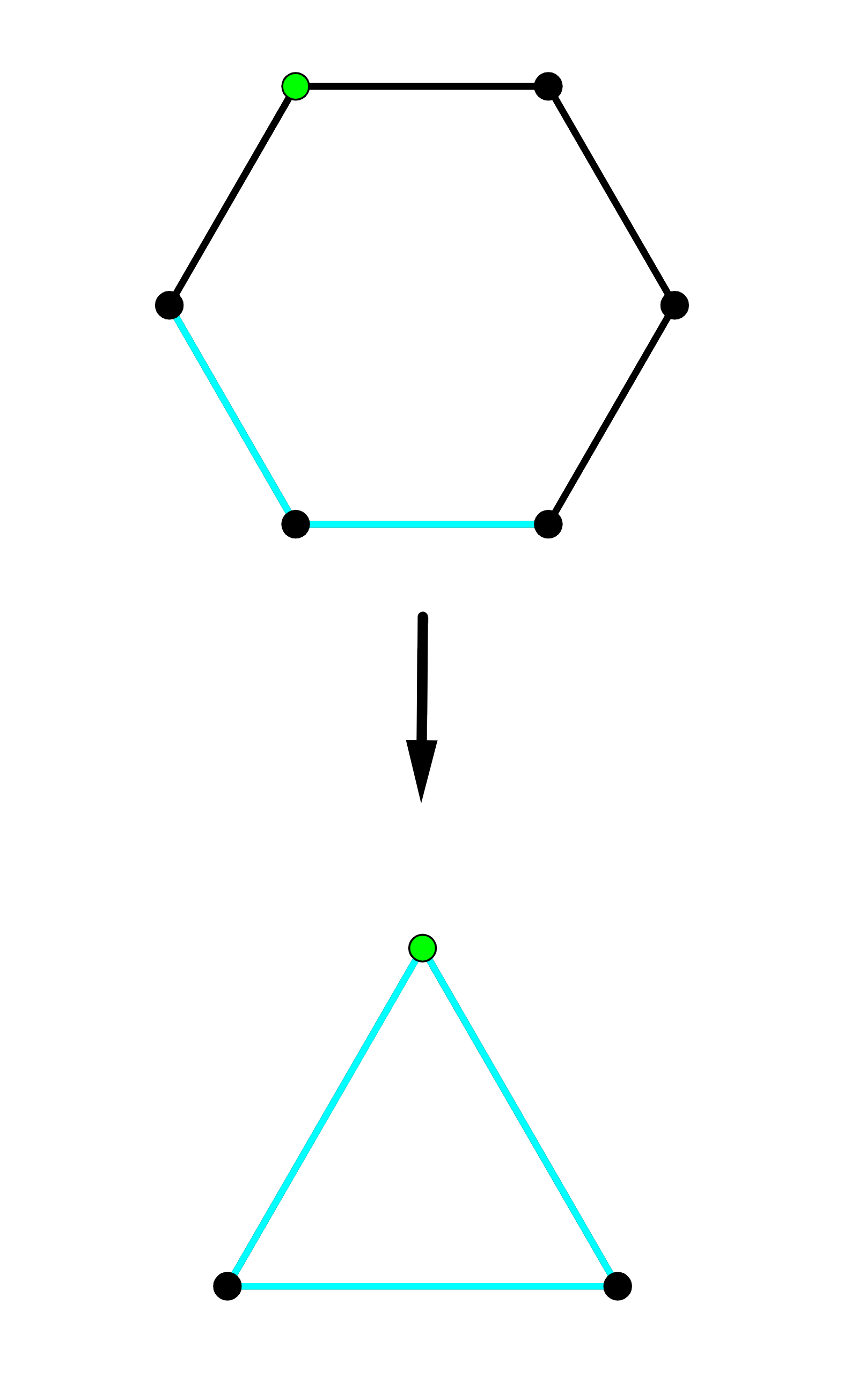}}
	
	\put(2.48,3.02){{\footnotesize $[0]$}}
	\put(3.25,3.02){{\footnotesize $[1]$}}
	\put(3.55,2.42){{\footnotesize $[2]$}}
	\put(3.25,1.82){{\footnotesize $[3]$}}
	\put(2.48,1.82){{\footnotesize $[4]$}}
	\put(2.19,2.42){{\footnotesize $[5]$}}

	\put(2.87,1.1){{\footnotesize $[0]$}}
	\put(2.33,0.2){{\footnotesize $[1]$}}
	\put(3.42,0.2){{\footnotesize $[2]$}}
	
	\end{picture}
	\caption{The local isomorphism $p: \C_{6} \to \C_{3}$}
	\label{Local_Iso}
\end{figure}
\end{ex}

\begin{df}\cite[p. 115]{GodsilAlgebraic}\label{Covering_Graph_Def}
Let $X$ and $\tX$ be graphs, and let $\rho: \tX \to X$ be a graph homomorphism. The pair $(\tX, \rho)$ is a \textit{covering graph of $X$} if $\rho$ is a surjective local isomorphism.
\end{df}

The pair $(C_6, \rho)$ in Example \ref{Local_Isomorphism_Ex} is also a covering graph of $C_3$. We often write the covering graph of a graph $X$ as simply $\rho: \tX \to X$.

\begin{df}\label{Lift_Def}
Let $X$ be a graph, and let $(\tX, \rho)$ be a covering graph of $X$. Given a graph homomorphism $f: Y \to X$, a \textit{lift of $f$} is a graph homomorphism $\widetilde{f}: Y \to \tX$ such that $\rho \circ \widetilde{f} = f$.
\end{df}

\begin{thm}[Path Lifting Property]\label{Path_Lifting_Property}
Let $(\tX, \rho)$ be a covering graph of $(X, x_0)$. For each stable graph homomorphism $f: I_\infty \to X$ with $f(n_f^1) = x_0$ and each vertex $\widetilde{x}_0 \in \rho^{-1}(x_0)$, there exists a unique lift $\widetilde{f}$ of $f$ starting at the vertex $\widetilde{x}_0$.
\[\begin{tikzcd}
& \tX \arrow[d, "\rho"] \\
I_\infty \arrow[r, "f"'] \arrow[ru, "\exists ! \widetilde{f}"] & X 
\end{tikzcd}\]
\end{thm}

\begin{proof}
Let $f$ and $\widetilde{x}_0$ be as in the statement. Define $\widetilde{f}: I_\infty \to \tX$ by $\widetilde{f}(i) = \widetilde{x}_0$ for all $i \leq n_f^1$ and recursively by
$$\widetilde{f}(i) = (\rho|_{N_{\widetilde{f}(i-1)}})^{-1}(f(i)) \;\;\;\;\; \text{for} \;\;\;\;\; i > n_f^1.$$ 
This means that $\widetilde{f}$ is defined using a different restriction for each $i \geq n_f^1$. It is not immediately obvious that this produces a graph homomorphism. Since $\widetilde{f}$ is defined recursively, we will use induction for $i \geq n_f^1$. 

By Lemma \ref{Relative_Local_Isomorphism_Lemma}, the graph homomorphism $(\rho|_{N_{\widetilde{f}(i-1)}})^{-1}$ exists. Since $f(i)$ is in the domain of $(\rho|_{N_{\widetilde{f}(i-1)}})^{-1}: N_{\rho(\widetilde{f}(i-1))} \to N_{\widetilde{f}(i-1)}$ for each $i \geq n_f^1$, $\widetilde{f}$ is well-defined.

For the base case, consider the edge $f(n_f^1) f(n_f^1 + 1) \in E(X)$. There is an edge $\widetilde{f}(n_f^1) \widetilde{f}(n_f^1 + 1) \in E(\tX)$ since $(\rho|_{N_{\widetilde{x}_0}})^{-1}: N_{x_0} \to N_{\widetilde{x}_0}$ is a graph homomorphism.
For the inductive hypothesis, suppose that $\widetilde{f}(i - 1) \widetilde{f}(i) \in E(\tX)$ for some $i > n_f^1$. By definition,
$$\widetilde{f}(i) = (\rho|_{N_{\widetilde{f}(i-1)}})^{-1}(f(i)) \;\;\;\;\; \text{and} \;\;\;\;\; \widetilde{f}(i + 1) = (\rho|_{N_{\widetilde{f}(i)}})^{-1}(f(i + 1)).$$
Here we have two different restrictions of $\rho$, however, $\widetilde{f}(i) \in N[\widetilde{f}(i - 1)] \cap N[\widetilde{f}(i)]$
by the inductive hypothesis. Thus
$$(\rho|_{N_{\widetilde{f}(i-1)}})^{-1}(f(i)) = (\rho|_{N_{\widetilde{f}(i)}})^{-1}(f(i)),$$
which implies that $\widetilde{f}$ is a graph homomorphism.

The map $\widetilde{f}$ is a lift of $f$, since by definition of $\widetilde{f}$ for all $i > n_f^1$,
$$\rho(\widetilde{f}(i)) = \rho((\rho|_{N_{\widetilde{f}(i-1)}})^{-1}(f(i))) = f(i).$$
To show that $\widetilde{f}$ is unique for each choice of $\widetilde{x}_0 \in \rho^{-1}(x_0)$, let $\widetilde{g}: I_\infty \to \tX$ be another graph homomorphism such that $\widetilde{g}(n_{\widetilde{g}}^1) = \widetilde{x}_0$ and $\rho \circ \widetilde{g} = f$. Thus $\rho(\widetilde{g}(i)) = f(i) = x_0$ for all $i \leq n_f^1$. By applying $(\rho|_{N_{\widetilde{f}(i)}})^{-1}$, we have that $\widetilde{g}(i) = \widetilde{x}_0$ for all $i \leq n_f^1$.
For the inductive hypothesis, suppose $\widetilde{g}(i) = \widetilde{f}(i)$ for some $i \geq n_f^1$. Since $\widetilde{f}$ and $\widetilde{g}$ are graph homomorphisms, this implies that $\widetilde{f}(i + 1), \widetilde{g}(i + 1) \in N[\widetilde{f}(i)]$. Since $\rho \circ \widetilde{g} = f = \rho \circ \widetilde{f}$, it follows that
$$\rho|_{N_{\widetilde{f}(i)}}(\widetilde{g}(i + 1)) = \rho|_{N_{\widetilde{f}(i)}}
(\widetilde{f}(i + 1)).$$
Thus by applying $(\rho|_{N_{\widetilde{f}(i)}})^{-1}$, we have that $\widetilde{g}(i + 1) = \widetilde{f}(i + 1)$.
\end{proof}

The Path Lifting Property (Theorem \ref{Path_Lifting_Property}) is used to prove the Homotopy Lifting Property (Theorem \ref{Homotopy_Lifting_Property}). However, the Homotopy Lifting Property does not hold for graphs with 3-cycles and 4-cycles. As mentioned previously, the 3-cycle and 4-cycle are A-contractible, but the cycles on five or more vertices are not. The analogous property in classical homotopy theory has no comparable condition. The statement of the Homotopy Lifting Property can be summarized by the following diagram.

\[\begin{tikzcd}
& & & \tX \arrow{ddd}{\rho} \\
 & & & \\
 & & & \\
 Y \arrow[rrr, shift left=1ex, "g"] \arrow[rrr, shift right=1ex, "f"', "\sim"] \arrow[uuurrr, dashed, shift left=2ex, " \widetilde{g}"] \arrow[uuurrr, "\widetilde{f}"', "\sim"] & & & X
\end{tikzcd}\]

Here $(\tX, \rho)$ is a covering graph of $X$, the $\sim$ between $f$ and $g$ represents the graph homotopy $H$ from $f$ to $g$, and the $\sim$ between $\widetilde{f}$ and $\widetilde{g}$ represents a lift $\widetilde{H}$ of $H$, a graph homotopy between a lift $\widetilde{f}$ of $f$ and a lift $\widetilde{g}$ of $g$. Thus if a lift $\widetilde{H}$ of $H$ exists, then a lift $\widetilde{g}$ of $g$ exists as well.

\begin{thm}[Homotopy Lifting Property]\label{Homotopy_Lifting_Property}
Let $X$ be a graph with no 3-cycles or 4-cycles and $(\tX, \rho)$ be a covering graph of $X$. Given a homotopy $H: Y \square I_n \to X$ from $f$ to $g$ and a lift $\widetilde{f}: Y \to \tX$ of $f$, there exists a unique homotopy $\widetilde{H}: Y \square I_n \to \tX$  that lifts $H$ with $\widetilde{H}(y, 0) = \widetilde{f}(y)$ for all $y \in V(Y)$.
\end{thm}

\begin{proof}
Let $X$, $(\tX, \rho)$, $H$ and $\widetilde{f}$ be as stated in the theorem. This proof is a combinatorial adaptation of the standard topological proof which can be found in \cite{HatcherAlgebraic}. 
Let $y \in V(Y)$ be an arbitrary vertex. First, we use the lift of $f$ to construct a lift of $H|_{N_y \square \{0\}} = f|_{N_y}$. Define $\widetilde{H}|_{N_y \square \{0\}} = \widetilde{f}|_{N_y}$. We now extend this map to $N_{(y, 0)}$, which contains the additional vertex $(y, 1)$. Since $\rho$ is a covering map, the restriction 
$$\rho|_{N_{\widetilde{H}(y,0)}}: N_{\widetilde{H}(y,0)} \to N_{\rho(\widetilde{H}(y, 0))}$$
is a bijection. Since $\rho(\widetilde{H}(y, 0)) = f(y) = H(y, 0)$, the inverse
$$(\rho|_{N_{\widetilde{H}(y,0)}})^{-1}: N_{H(y,0)} \to N_{\widetilde{H}(y,0)}$$
exists by Lemma \ref{Relative_Local_Isomorphism_Lemma}. Since $H$ is a graph homomorphism, there is an inclusion of sets $H(N[y, 0]) \subseteq N[H(y, 0)]$. In particular, $H(y, 1) \in N[H(y, 0)]$, which implies that $H(y, 1)$ is in the domain of $(\rho|_{N_{\widetilde{H}(y,0)}})^{-1}$. Define $\widetilde{H}(y, 1) = (\rho|_{N_{\widetilde{H}(y,0)}})^{-1}(H(y, 1))$. Since $\widetilde{f}$ is a lift of $f$ and $H|_{N_y \square \{0\}} = f|_{N_y}$, it follows that $\widetilde{f}|_{N_y} = (\rho|_{N_{\widetilde{H}(y,0)}})^{-1} \circ  H|_{N_y \square \{0\}}$. Thus we have defined $\widetilde{H}|_{N(y,0)}$, and it is a graph homomorphism because it is the composition of graph homomorphisms.

We now proceed by recursion to define a lift of $H|_{N_{(y,i+1)}}$ for each $0 \leq i < n$, which agrees with the previous lift of $H|_{N_{(y,i)}}$ on $N[y,i-1] \cap N[y,i]$.
	\begin{figure}[ht]
	\begin{picture}(5,1.25)
	\put(0,0){\includegraphics[scale=0.105]{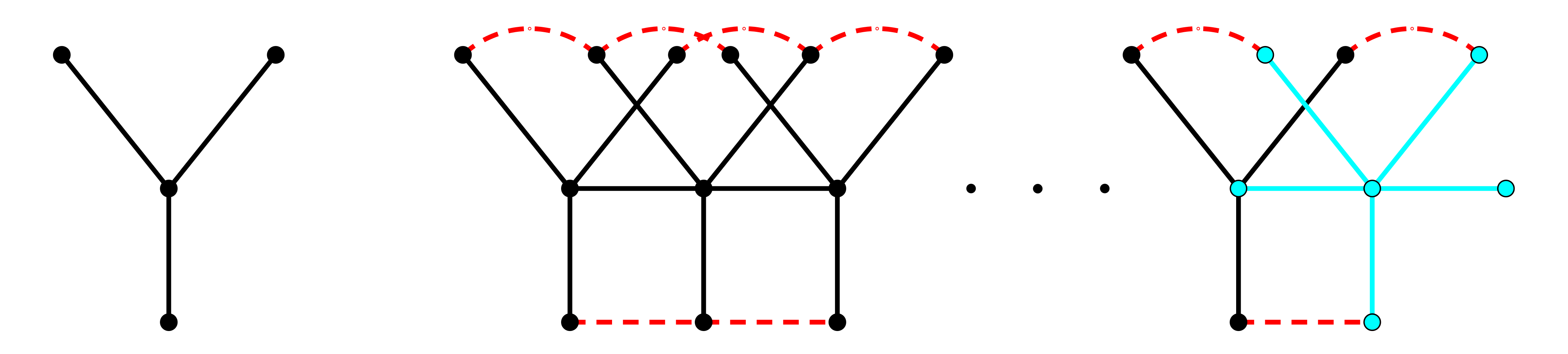}}
	\put(0.45,1.13){$N_y$}
	\put(0.56,0.45){{\small $y$}}

	\put(1.5,0.42){{\tiny $(y,0)$}}
	\put(1.91,0.42){{\tiny $(y,1)$}}
	\put(2.32,0.42){{\tiny $(y,2)$}}

	\put(3.4,0.42){{\tiny $(y,i-1)$}}
	\put(4,0.42){{\tiny $(y,i)$}}
	\put(4.45,0.42){{\tiny $(y,i+1)$}}
	\end{picture}
	\caption{Graph $N_{(y, 0)} \cup \cdots \cup N_{(y, i-1)} \cup N_{(y, i)}$}
	\label{Union_0_to_i_Fig}
	\end{figure}
Figure \ref{Union_0_to_i_Fig} illustrates the graph $N_{(y,0)} \cup \cdots \cup N_{(y,i)}$, in the case that the vertex $y$ has three adjacent vertices. The subgraph $N_{(y,i)}$ is shown in light blue, and the dashed edges shown in red are the edges in $N_y \square I_i$, which are not in $N_{(y,0)} \cup \cdots \cup N_{(y,i)}$. 

Assume that $H|_{N_{(y,0)} \cup \cdots \cup N_{(y,i)}}$ has a lift $\widetilde{H}|_{N_{(y,0)} \cup \cdots \cup N_{(y,i)}}$ for some $0 \leq i < n$. Since $(y, i + 1) \in N[y, i]$, it follows that $\widetilde{H}(y, i + 1)$ is defined. Define
$$\widetilde{H}|_{N_{(y,i+1)}} = (\rho|_{N_{\widetilde{H}(y,i+1)}})^{-1} \circ H|_{N_{(y,i+1)}}$$ 
for each $0 \leq i < n$.
Showing that $\widetilde{H}|_{N_{(y,i+1)}}$ is well-defined and a graph homomorphism is very similar to the base case
and is illustrated by the following diagram.
	\[ \begin{tikzcd}
	& & & N_{\widetilde{H}(y, i+1)} \\
	& & & \\
	& & & \\
	N_{(y, i+1)} \arrow[uuurrr, "\widetilde{H}|_{N_{(y, i+1)}}"] \arrow[rr, "H|_{N_{(y,i+1)}}"'] & & H(N_{(y, i+1)}) \arrow[r, hook] & N_{H(y, i+1)} \arrow[uuu, "(\rho|_{N_{\widetilde{H}(y, i+1)}})^{-1}"']
	\end{tikzcd} \]
Thus after a finite number of steps, the lift $\widetilde{H}|_{N_{(y,0)} \cup \cdots \cup N_{(y,n)}}$ is defined.

For all $x \in N[y]$, there is an edge $(x, i) (x, i + 1) \in E(N_y \square I_n)$ for $0 \leq i < n$. In order to extend $\widetilde{H}|_{N_{(y,0)} \cup \cdots \cup N_{(y,n)}}$ to a graph homomorphism with domain $N_y \square I_n$, there must be an edge $\widetilde{H}(x, i) \widetilde{H}(x, i+ 1) \in E(\tX)$ for all $0 \leq i < n$. By definition of $\widetilde{H}|_{N_{(y,0)} \cup \cdots \cup N_{(y,n)}}$,
$$\widetilde{H}(x, i) = (\rho|_{N_{\widetilde{H}}(y,i)})^{-1} \circ H|_{N(y,i)}(x, i) = (\rho|_{N_{\widetilde{H}(y,i)}})^{-1} \circ H(x, i)$$
and
$$\widetilde{H}(x, i + 1) = (\rho|_{N_{\widetilde{H}(y,i+1)}})^{-1} \circ H|_{N(y,i+1)}(x, i + 1) = (\rho|_{N_{\widetilde{H}(y,i+1)}})^{-1} \circ H(x, i + 1).$$
In order to show that there is an edge $\widetilde{H}(x, i) \widetilde{H}(x, i+1) \in E(\tX)$ for all $0 \leq i < n$, we will examine the 4-cycle of $N_y \square I_n$ shown in light blue in Figure \ref{Union_0_to_n_Fig}.
    \begin{figure}[ht]
	\begin{picture}(6,1.2)
	\put(0.06,0){\includegraphics[scale=0.265]{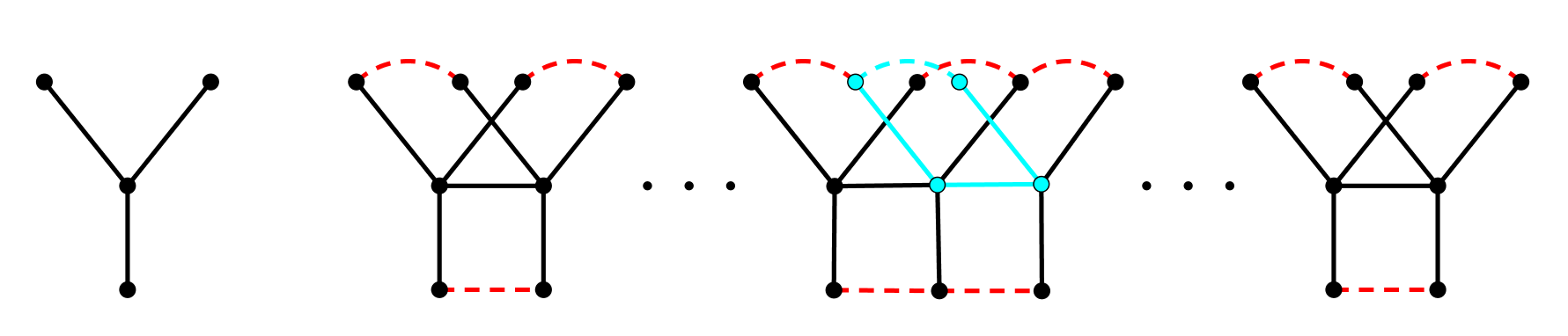}}
	\put(0.4,0.9){$N_{y}$}
	\put(0.32,0.45){{\small $y$}}
	\put(0.05,0.78){{\small $x$}}

	\put(1.18,0.38){{\tiny $(y,0)$}}
	\put(1.51,0.38){{\tiny $(y,1)$}}

	\put(2.25,0.38){{\tiny $(y,i-1)$}}
	\put(2.758,0.38){{\tiny $(y,i)$}}
	\put(3.33,0.38){{\tiny $(y,i+1)$}}

	\put(2.55,0.9){{\tiny $(x,i)$}}
	\put(2.95,0.9){{\tiny $(x,i+1)$}}

	\put(3.8,0.38){{\tiny $(y,n-1)$}}
	\put(4.3,0.38){{\tiny $(y,n)$}}
	\end{picture}
	\caption{Union of neighborhoods $N_{(y, 0)} \cup \cdots \cup N_{(y, n)}$}
	\label{Union_0_to_n_Fig}
	\end{figure}
We denote this 4-cycle subgraph by $C_{x,i}$. Since $H$ is a graph homomorphism and $X$ contains no 3-cycles or 4-cycles, there are only nine ways that $H$ maps $C_{x,i}$ to $X$, illustrated in Figure \ref{Cases_Fig}. The label `=’ means that H maps the pair of vertices to the same vertex in $X$. The label `a' means that $H$ maps the pair of vertices to distinct adjacent vertices in $X$. In cases (8) and (9), the pair of vertices being mapped to the same vertex are circled in red.
	\begin{figure}[ht]
	\begin{picture}(4, 4.5)
	\put(-0.35,0){\includegraphics[scale=1.1]{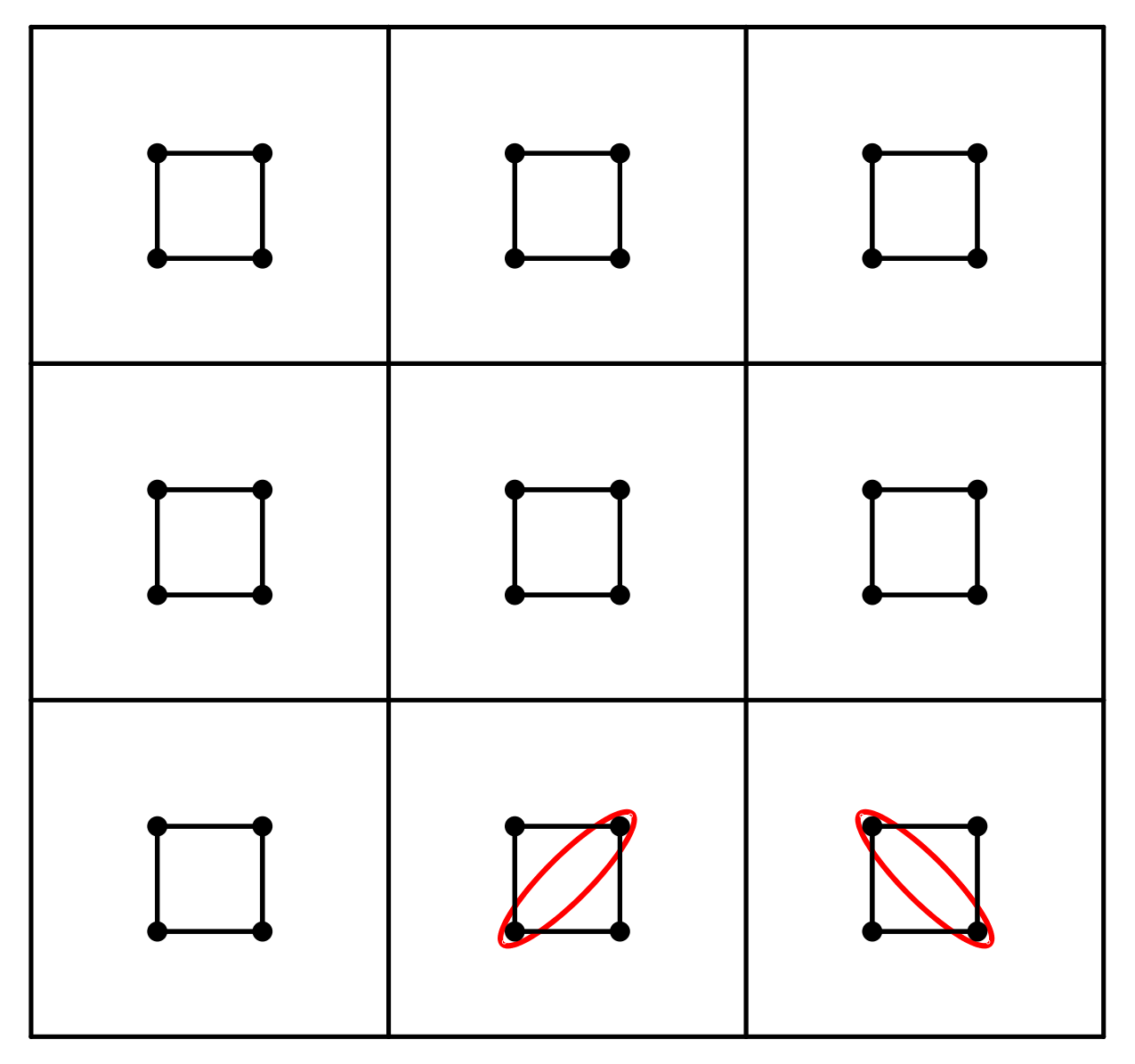}}
	\put(-0.18,4.12){{\small $(1)$}}
	\put(1.3,4.12){{\small $(2)$}}
	\put(2.78,4.12){{\small $(3)$}}
	\put(-0.18,2.73){{\small $(4)$}}
	\put(1.3,2.73){{\small $(5)$}}
	\put(2.78,2.73){{\small $(6)$}}
	\put(-0.18,1.34){{\small $(7)$}}
	\put(1.3,1.34){{\small $(8)$}}
	\put(2.78,1.34){{\small $(9)$}}

	\put(0.03,3.81){{\tiny $(x, i)$}}
	\put(0.03,3.2){{\tiny $(y, i)$}}
	\put(0.78,3.81){{\tiny $(x, i+1)$}}
	\put(0.78,3.2){{\tiny $(y, i+1)$}}

	\put(1.505,3.81){{\tiny $(x, i)$}}
	\put(1.505,3.2){{\tiny $(y, i)$}}
	\put(2.255,3.81){{\tiny $(x, i+1)$}}
	\put(2.255,3.2){{\tiny $(y, i+1)$}}

	\put(2.98,3.81){{\tiny $(x, i)$}}
	\put(2.98,3.2){{\tiny $(y, i)$}}
	\put(3.73,3.81){{\tiny $(x, i+1)$}}
	\put(3.73,3.2){{\tiny $(y, i+1)$}}

	\put(0.03,2.42){{\tiny $(x, i)$}}
	\put(0.03,1.81){{\tiny $(y, i)$}}
	\put(0.78,2.42){{\tiny $(x, i+1)$}}
	\put(0.78,1.81){{\tiny $(y, i+1)$}}

	\put(1.505,2.42){{\tiny $(x, i)$}}
	\put(1.505,1.81){{\tiny $(y, i)$}}
	\put(2.255,2.42){{\tiny $(x, i+1)$}}
	\put(2.255,1.81){{\tiny $(y, i+1)$}}
	
	\put(2.98,2.42){{\tiny $(x, i)$}}
	\put(2.98,1.81){{\tiny $(y, i)$}}
	\put(3.73,2.42){{\tiny $(x, i+1)$}}
	\put(3.73,1.81){{\tiny $(y, i+1)$}}
	
	\put(0.03,1.03){{\tiny $(x, i)$}}
	\put(0.03,0.42){{\tiny $(y, i)$}}
	\put(0.78,1.03){{\tiny $(x, i+1)$}}
	\put(0.78,0.42){{\tiny $(y, i+1)$}}

	\put(1.505,1.03){{\tiny $(x, i)$}}
	\put(1.505,0.42){{\tiny $(y, i)$}}
	\put(2.255,1.05){{\tiny $(x, i+1)$}}
	\put(2.255,0.42){{\tiny $(y, i+1)$}}
	
	\put(2.98,1.03){{\tiny $(x, i)$}}
	\put(2.98,0.42){{\tiny $(y, i)$}}
	\put(3.73,1.03){{\tiny $(x, i+1)$}}
	\put(3.73,0.42){{\tiny $(y, i+1)$}}

	\put(0.14,3.505){{\tiny $=$}}
	\put(0.47,3.81){{\tiny $=$}}
	\put(0.82,3.505){{\tiny $=$}}
	\put(0.47,3.2){{\tiny $=$}}

	\put(1.6,3.505){{\tiny $=$}}
	\put(1.935,3.81){{\tiny a}}
	\put(2.3,3.505){{\tiny a}}
	\put(1.94,3.2){{\tiny $=$}}

	\put(3.06,3.505){{\tiny $=$}}
	\put(3.415,3.81){{\tiny a}}
	\put(3.78,3.505){{\tiny $=$}}
	\put(3.415,3.2){{\tiny a}}
	
	\put(0.14,2.115){{\tiny $=$}}
	\put(0.47,2.42){{\tiny $=$}}
	\put(0.82,2.115){{\tiny a}}
	\put(0.475,1.81){{\tiny a}}
	
	\put(1.605,2.115){{\tiny a}}
	\put(1.94,2.42){{\tiny a}}
	\put(2.3,2.115){{\tiny $=$}}
	\put(1.94,1.81){{\tiny $=$}}
	
	\put(3.06,2.115){{\tiny a}}
	\put(3.415,2.42){{\tiny $=$}}
	\put(3.78,2.115){{\tiny a}}
	\put(3.415,1.81){{\tiny $=$}}
	
	\put(0.14,0.725){{\tiny a}}
	\put(0.47,1.03){{\tiny $=$}}
	\put(0.82,0.725){{\tiny $=$}}
	\put(0.47,0.42){{\tiny a}}
	
	\put(1.605,0.725){{\tiny a}}
	\put(1.945,1.03){{\tiny a}}
	\put(2.3,0.725){{\tiny a}}
	\put(1.945,0.42){{\tiny a}}
	\put(1.945,0.7){\rotatebox{45}{\scriptsize $=$}}
	
	\put(3.065,0.725){{\tiny a}}
	\put(3.415,1.03){{\tiny a}}
	\put(3.78,0.725){{\tiny a}}
	\put(3.415,0.42){{\tiny a}}
	\put(3.42,0.745){\rotatebox{135}{\scriptsize $=$}}
	\end{picture}
	\caption{The cases of how $H$ maps $C_{x, i}$ to $X$}
	\label{Cases_Fig}
	\end{figure}

For cases (1)-(8), there is an inclusion of sets $V(H(C_{x,i})) \subseteq N[H(y, i)]$, and there is an edge $H(x, i) H(x, i+1) \in E(X)$. Thus the subgraph $C_{x,i}$ is mapped by $H$ into the domain of the inverse $(\rho|_{N_{\widetilde{H}(y,i)}})^{-1}: N_{H(y,i)} \to N_{\widetilde{H}(y,i)}$. Since we have the composition $\widetilde{H}|_{C_{x,i}} = (\rho|_{N_{\widetilde{H}(y,i)}})^{-1} \circ H|_{C_{x,i}}$, it follows that there is an edge $\widetilde{H}(x, i) \widetilde{H}(x, i + 1) \in E(\tX)$.

For case (9), there is an inclusion of sets $H(C_{x,i}) \subset N[H(y, i + 1)]$, and there is an edge $H(x, i) H(x, i + 1) \in E(X)$. Thus the subgraph $C_{x,i}$ is mapped by $H$ into the domain of the inverse $(\rho|_{N_{\widetilde{H}(y,i+1)}})^{-1}: N_{H(y,i+1)} \to N_{\widetilde{H}(y,i+1)}$. Since we have the composition $\widetilde{H}|_{C_{x,i}} = (\rho|_{N_{\widetilde{H}(y,i+1)}})^{-1} \circ H|_{C_{x,i}}$, it follows that there is an edge $\widetilde{H}(x, i) \widetilde{H}(x, i + 1) \in E(\tX)$. Thus we can extend the graph homomorphism $\widetilde{H}|_{N(y,0) \cup \cdots \cup N(y,n)}$ to $\widetilde{H}|_{N_y \square I_n}$.

The restriction $H|_{\{y\} \square I_n}$ is a graph homomorphism from $I_n$ to $X$. By the Uniqueness of Path Lifting (\ref{Path_Lifting_Property}), the lift $\widetilde{H}|_{\{y\} \square I_n}: I_n \to \tX$ of $H|_{\{y\} \square I_n}$ is unique with $\widetilde{H}|_{\{y\} \square I_n}(0) = \widetilde{H}(y, 0) = \widetilde{f}(y)$. Since each graph homomorphism $H|_{\{x\} \square I_n}: I_n \to X$ must have a unique lift $\widetilde{H}|_{\{x\} \square I_n}: I_n \to \tX$ for all $x \in N[y]$ with $\widetilde{H}|_{\{x\} \square I_n}(0) = \widetilde{H}(x, 0) = \widetilde{f}(x)$, the lift $\widetilde{H}|_{N_{y \square I_n}}$ must be unique for each $y \in V(Y)$. Since $\widetilde{H}_x$ is unique for each $x \in V(Y)$ and is a restriction of the graph homomorphism $\widetilde{H}|_{N_{y \square I_n}}$ for each $y \in V(Y)$ such that $x \in N[y]$, the graph homomorphisms $\widetilde{H}|_{N_y \square I_n}$ must form a unique lift $\widetilde{H}$ of the homotopy $H$.
\end{proof}

Recall that the Homotopy Lifting Property (Theorem \ref{Homotopy_Lifting_Property}) does not hold for graphs containing 3-cycles and 4-cycles. The following example demonstrates how Theorem \ref{Homotopy_Lifting_Property} fails for the cycle $C_4$.

\begin{ex}\label{Homotopy_Lift_Fail_Ex}
Let $f: I_4 \to \C_4$ be the graph homomorphism defined by 
$$f(0) = [0], \; f(1) = [1], \; f(2) = [2], \; f(3) = [3], \;\; \textrm{and} \;\; f(4) = [0].$$
Let $g: I_4 \to \C_4$ be the graph homomorphism defined by $g(i) = [0]$ for all $i \in \{0, 1, 2, 3, 4\}$. The pair $(I_\infty, \rho_4)$ is a covering graph of $\C_4$, where $\rho_4$ is defined by $\rho_4(i) = [i \mod 4]$ for $i \in \Z$. Figure \ref{Homotopy_C4_Fig} depicts a graph homotopy $H: I_4 \square I_2 \to \C_4$ from $f$ to $g$.
	\begin{figure}[h!]
	\begin{picture}(4.3, 2.4)
	\put(0,-0.2){\includegraphics[scale=1.27]{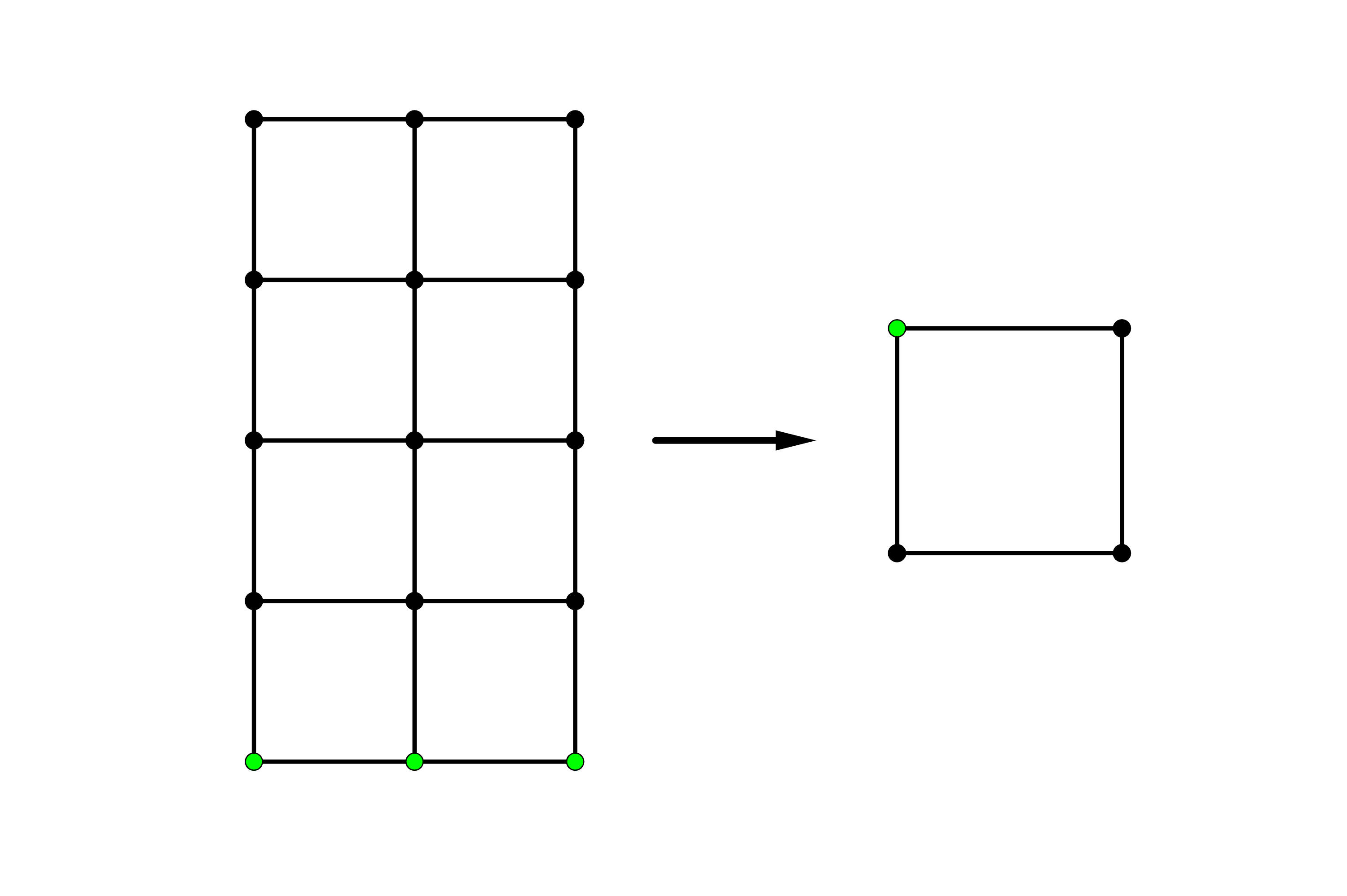}}

	\put(0.6,2.17){{\small $4$}}
	\put(0.6,1.68){{\small $3$}}
	\put(0.6,1.18){{\small $2$}}
	\put(0.6,0.69){{\small $1$}}
	\put(0.6,0.18){{\small $0$}}

	\put(0.74,2.33){{\small $0$}}
	\put(1.25,2.33){{\small $1$}}
	\put(1.76,2.33){{\small $2$}}

	\put(2.19,1.3){$H$}

	\put(2.61,1.62){{\small $[0]$}}
	\put(3.54,1.62){{\small $[1]$}}
	\put(3.54,0.75){{\small $[2]$}}
	\put(2.61,0.75){{\small $[3]$}}

	\put(0.82,2.11){\color{red}{\scriptsize $[0]$}}
	\put(0.82,1.61){\color{red}{\scriptsize $[3]$}}
	\put(0.82,1.11){\color{red}{\scriptsize $[2]$}}
	\put(0.82,0.61){\color{red}{\scriptsize $[1]$}}
	\put(0.82,0.11){\color{red}{\scriptsize $[0]$}}

	\put(1.32,2.11){\color{red}{\scriptsize $[0]$}}
	\put(1.32,1.61){\color{red}{\scriptsize $[0]$}}
	\put(1.32,1.11){\color{red}{\scriptsize $[1]$}}
	\put(1.32,0.61){\color{red}{\scriptsize $[1]$}}
	\put(1.32,0.11){\color{red}{\scriptsize $[0]$}}

	\put(1.82,2.11){\color{red}{\scriptsize $[0]$}}
	\put(1.82,1.61){\color{red}{\scriptsize $[0]$}}
	\put(1.82,1.11){\color{red}{\scriptsize $[0]$}}
	\put(1.82,0.61){\color{red}{\scriptsize $[0]$}}
	\put(1.82,0.11){\color{red}{\scriptsize $[0]$}}
	\end{picture}
	\caption{A homotopy $H:I_{4} \square I_{2} \to \C_{4}$}
	\label{Homotopy_C4_Fig}
	\end{figure}
	
By the Path Lifting Property (Theorem \ref{Path_Lifting_Property}), there is a unique lift $\widetilde{f}: I_4 \to I_\infty$ of $f$ given by 
$$\widetilde{f}(0) = 0, \; \widetilde{f}(1) = 1, \; \widetilde{f}(2) = 2, \; \widetilde{f}(3) = 3, \;\; \textrm{and} \;\; \widetilde{f}(4) = 4.$$
We can attempt to construct a lift of $H$ using the path lifts of the restrictions $H|_{\{i\} \square I_2}: I_2 \to C_4$ for each $0 \leq i \leq 4$. The constructed map $\widetilde{H}: I_4 \square I_2 \to I_\infty$ is depicted in Figure \ref{Homotopy_C4_Lift_Fig} and is not a graph homomorphism. The edges shown in red are incident to vertices that are not mapped to adjacent vertices of $I_\infty$.
	\begin{figure}[ht]
	\begin{picture}(4, 2.82)
	\put(0,-0.2){\includegraphics[scale=1.25]{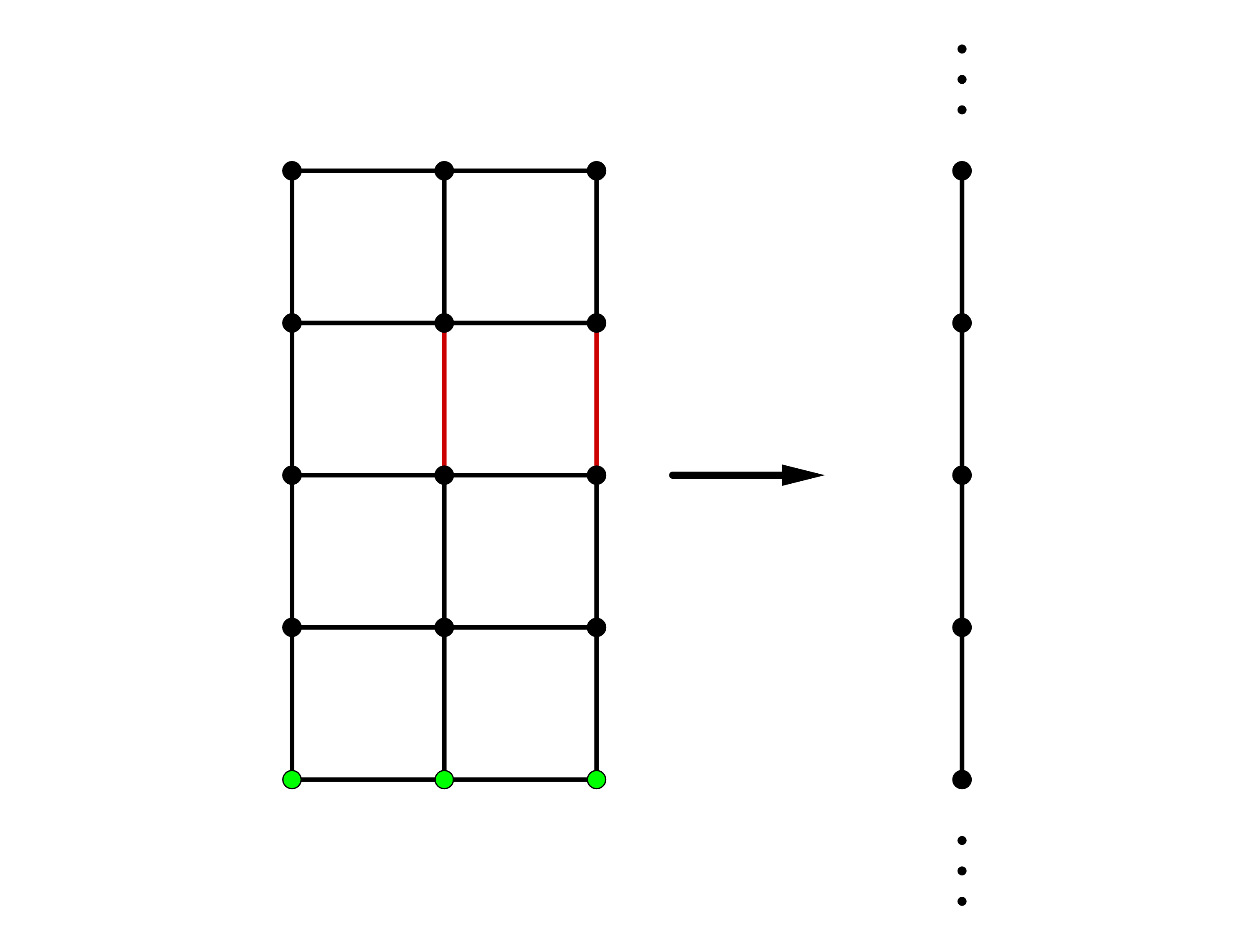}}
	\put(0.8,2.3){{\small $4$}}
	\put(0.8,1.8){{\small $3$}}
	\put(0.8,1.31){{\small $2$}}
	\put(0.8,0.82){{\small $1$}}
	\put(0.8,0.32){{\small $0$}}

	\put(0.9,2.45){{\small $0$}}
	\put(1.4,2.45){{\small $1$}}
	\put(1.9,2.45){{\small $2$}}

	\put(2.3,1.4){$\widetilde{H}$}

	\put(2.95,2.3){{\small $4$}}
	\put(2.95,1.8){{\small $3$}}
	\put(2.95,1.31){{\small $2$}}
	\put(2.95,0.82){{\small $1$}}
	\put(2.95,0.32){{\small $0$}}

	\put(0.98,2.2){\color{red}{\footnotesize $4$}}
	\put(0.98,1.7){\color{red}{\footnotesize $3$}}
	\put(0.98,1.21){\color{red}{\footnotesize $2$}}
	\put(0.98,0.72){\color{red}{\footnotesize $1$}}
	\put(0.98,0.22){\color{red}{\footnotesize $0$}}

	\put(1.48,2.2){\color{red}{\footnotesize $4$}}
	\put(1.48,1.7){\color{red}{\footnotesize $4$}}
	\put(1.48,1.21){\color{red}{\footnotesize $1$}}
	\put(1.48,0.72){\color{red}{\footnotesize $1$}}
	\put(1.48,0.22){\color{red}{\footnotesize $0$}}

	\put(1.98,2.2){\color{red}{\footnotesize $4$}}
	\put(1.98,1.7){\color{red}{\footnotesize $4$}}
	\put(1.98,1.21){\color{red}{\footnotesize $0$}}
	\put(1.98,0.72){\color{red}{\footnotesize $0$}}
	\put(1.98,0.22){\color{red}{\footnotesize $0$}}
	\end{picture}
	\caption{The Map $\widetilde{H}: I_{4} \square I_{2} \to I_{\infty}$}
	\label{Homotopy_C4_Lift_Fig}
	\end{figure}

Since $H$ is a graph homomorphism from $I_4 \square I_2$ to $C_4$, we could also try to construct a lift of $H$ with the lifts of the restrictions $H|_{I_4 \square \{j\}}: I_4 \to C_4$. However, this construction also fails to be a graph homomorphism.
\end{ex}

We now use the Path Lifting Property (Theorem \ref{Path_Lifting_Property}) and the Homotopy Lifting Property (Theorem \ref{Homotopy_Lifting_Property}) to prove the general Lifting Criterion (Theorem \ref{Lifting_Criterion}).

\begin{df}\label{Connected_Def}
A graph $X$ is \textit{connected} if for each $v, w \in V(X)$, there exists a stable graph homomorphism $f \in S_1(X)$ such that $f(n_f^1) = v$ and $f(p_f^1) = w$.
\end{df}

This is an adaptation of the standard definition of a connected graph found in e.g. \cite{WestGraph}, to involve a stable graph homomorphism. We will use
Lemma \ref{Induced_Map_Composition_Lemma} and Definition \ref{Connected_Def} as well as the Path Lifting Property (Theorem \ref{Path_Lifting_Property})
and Homotopy Lifting Property (Theorem \ref{Homotopy_Lifting_Property}) in the proof of this last result, the Lifting
Criterion (Theorem \ref{Lifting_Criterion}).

\begin{thm}[Lifting Criterion]\label{Lifting_Criterion}
Let $(\tX, \rho)$ be a covering graph of a graph $X$, and let $f: (Y, y_0) \to (X, x_0)$ be a graph homomorphism with a connected graph $Y$. If $X$ contains neither 3-cycles nor 4-cycles, then there is a lift $\widetilde{f}: (Y, y_0) \to (\tX, \widetilde{x}_0)$ of $f$ if and only if $f_\ast(A_1(Y, y_0)) \subseteq \rho_\ast(A_1(\tX, \widetilde{x}_0))$.
\end{thm}

\[\begin{tikzcd}
& & (\tX, \widetilde{x}_0) \arrow{dd}{\rho} \\
& & \\
(Y, y_0) \arrow[uurr, "\widetilde{f}"]  \arrow[rr, "f"'] & & (X, x_0)
\end{tikzcd}\]

\begin{proof}
Let $X$, $(\tX, \rho)$, and $f$ be as in the statement of the theorem. Suppose a lift $\widetilde{f}: (Y, y_0) \to (\tX, \widetilde{x}_0)$ of $f$ exists. Then $\rho \circ \widetilde{f} = f$, which implies that $\rho_\ast \circ \widetilde{f}_\ast = f_\ast$ by Lemma \ref{Induced_Map_Composition_Lemma}. 
It follows immediately that
$$f_\ast(A_1(Y, y_0)) \subseteq \rho_\ast(A_1(\tX, \widetilde{x}_0)).$$

Conversely, suppose $f_\ast(A_1(Y, y_0)) \subseteq \rho_\ast(A_1(\tX, \widetilde{x}_0))$. Let $y \in V(Y)$. Since $Y$ is connected, there is a stable graph homomorphism $\gamma_y: I_\infty \to Y$ with $\gamma_y(n_{\gamma_y}^1) = y_0$ and $\gamma_y(p_{\gamma_y}^1) = y$. Thus $f \circ \gamma_y: I_\infty \to X$ is a stable graph homomorphism with $f(\gamma_y(n_{\gamma_y}^1)) = x_0$ and $f(\gamma_y(p_{\gamma_y}^1)) = f(y) \in V(X)$. We write $f \circ \gamma_y$ as $f \gamma_y$ for the remainder of this proof to avoid cumbersome notation. By the Path Lifting Property (\ref{Path_Lifting_Property}), there is a unique lift $\widetilde{f \gamma_y} : I_\infty \to \tX$ with $\widetilde{f \gamma_y}(n_{\gamma_y}^1) = \widetilde{x}_0$. Define $\widetilde{f}: Y \to \tX$ by $\widetilde{f}(y) = \widetilde{f \gamma_y}(p_{\gamma_y}^1) \in \rho^{-1}(f(y))$.
\[\begin{tikzcd}
& & & (\tX, \widetilde{x}_0) \arrow{dd}{\rho} \\
& & & \\
I_\infty \arrow[r, "\gamma_y"'] \arrow[uurrr, bend left, "\widetilde{f \gamma_y}"] & (Y, y_0) \arrow[uurr, dashed, "\exists \widetilde{f}"]  \arrow[rr,  "f"'] & & (X, x_0)
\end{tikzcd}\]

Note that the choice of $\gamma_y$ is not unique. We must prove that $\widetilde{f}(y)$ does not depend on the choice of $\gamma_y$ to show that $\widetilde{f}$ is well-defined.

Suppose $\beta_y: I_\infty \to Y$ is another stable graph homomorphism with $\beta_y(n_{\beta_y}^1) = y_0$ and $\beta_y(p_{\beta_y}^1) = y$. Then $f \circ \beta_y: I_\infty \to X$ is a stable graph homomorphism with $f(\beta_y(n_{\beta_y}^1)) = x_0$ and $f(\beta_y(p_{\beta_y}^1)) = f(y)$. Similarly, we write $f \circ \beta_y$ as $f \beta_y$ for the remainder of the proof. In this proof, we will use the Path Lifting Property (\ref{Path_Lifting_Property}) and the Homotopy Lifting Property (\ref{Homotopy_Lifting_Property}) to show that $\widetilde{f \gamma_y}(p_{\gamma_y}^1) = \widetilde{f \beta_y}(p_{\beta_y}^1)$.

Recall that $\overline{\gamma_y}: I_\infty \to Y$ is defined by $\overline{\gamma_y}(i) = \gamma_y(-i)$ for all $i \in \Z$. Since $\overline{\gamma_y}(n_{\overline{\gamma_y}}^1) = y = \beta_y(p_{\beta_y}^1)$, the concatenation $\overline{\gamma_y} \cdot \beta_y: I_\infty \to Y$ is a graph homomorphism. By Lemma \ref{Concatenation_Prop}, $\overline{\gamma_y} \cdot \beta_y$ stabilizes in the negative direction at $n_{\overline{\gamma_y} \cdot \beta_y}^1 = n_{\beta_y}^1 - p_{\beta_y}^1$ and in the positive direction at $p_{\overline{\gamma_y} \cdot \beta_y}^1 = p_{\overline{\gamma_y}}^1 - n_{\overline{\gamma_y}}^1 = -n_{\gamma_y}^1 + p_{\gamma_y}^1$. With the definitions of concatenation and $\overline{\gamma}_y$, it can be verified that
\begin{eqnarray*}
\overline{\gamma_y} \cdot \beta_y(i) 
& = & \begin{cases}
\gamma_y(-i + p_{\gamma_y}^1) & \text{for} \;\;\; i \geq 0, \\
\beta_y(i + p_{\beta_y}^1) & \text{for} \;\;\; i \leq 0,
\end{cases} 
\end{eqnarray*}
Therefore,
\begin{eqnarray*}
\overline{\gamma_y} \cdot \beta_y(n_{\overline{\gamma_y} \cdot \beta_y}^1) 
& = & \overline{\gamma_y} \cdot \beta_y(n_{\beta_y}^1 - p_{\beta_y}^1) \\
& = & \beta_y(n_{\beta_y}^1 - p_{\beta_y}^1 + p_{\beta_y}^1) \\
& = & \beta_y(n_{\beta_y}^1) \\
& = & y_0
\end{eqnarray*}
and
\begin{eqnarray*}
\overline{\gamma_y} \cdot \beta_y(p_{\overline{\gamma_y} \cdot \beta_y}^1) 
& = & \overline{\gamma_y} \cdot \beta_y(-n_{\gamma_y}^1 + p_{\gamma_y}^1) \\
& = & \gamma_y(n_{\gamma_y}^1 - p_{\gamma_y}^1 + p_{\gamma_y}^1) \\
& = & \gamma_y(n_{\gamma_y}^1) \\
& = & y_0.
\end{eqnarray*}
Thus $[\overline{\gamma_y} \cdot \beta_y] \in A_1(Y, y_0)$, namely, $\overline{\gamma_y} \cdot \beta_y$ is a `loop’ in the graph $Y$ based at the distinguished vertex $y_0$. Since $f_\ast$ is a group homomorphism by Lemma \ref{Induced_Map_Group_Homomorphism_Lemma},
$$f_\ast([\overline{\gamma_y} \cdot \beta_y]) = [f \circ (\overline{\gamma_y} \cdot \beta_y)] = [\overline{f \gamma_y} \cdot f \beta_y].$$ 
Therefore, $[\overline{f \gamma_y} \cdot f \beta_y] \in f_\ast(A_1(Y, y_0)) \subseteq \rho_\ast(A_1(\tX, \widetilde{x}_0))$, which implies that there exists
an equivalence class $[g] \in A_1(\tX, \widetilde{x}_0)$ such that $\rho_\ast([g]) = [\overline{f \gamma_y} \cdot f \beta_y]$. Thus there exists a graph homotopy from $\rho \circ g$ to $\overline{f \gamma_y} \cdot f \beta_y$. By the Path Lifting Property (\ref{Path_Lifting_Property}), there is a unique lift
$$\widetilde{\overline{f \gamma_y} \cdot f \beta_y}: I_\infty \to \tX$$
of $\overline{f \gamma_y} \cdot f \beta_y$ with $\widetilde{\overline{f \gamma_y} \cdot f \beta_y}(n_{\widetilde{\overline{f \gamma_y} \cdot f \beta_y}}^1) = \widetilde{x}_0$. Since $X$ contains neither 3-cycles nor 4-cycles, the Homotopy Lifting Property (Theorem \ref{Homotopy_Lifting_Property}) holds. Thus there exists a lifted homotopy from $g$, which is a unique lift of $p \circ g$, to $\widetilde{\overline{f \gamma_y} \cdot f \beta_y}$. Since $[g] \in A_1(\tX, \widetilde{x}_0)$, it follows that
$$\widetilde{\overline{f \gamma_y} \cdot f \beta_y}(n_{\overline{f \gamma_y} \cdot f \beta_y}^1) = \widetilde{\overline{f \gamma_y} \cdot f \beta_y}(p_{\overline{f \gamma_y} \cdot f \beta_y}^1) = \widetilde{x}_0$$
as well. 
By the uniqueness of the Path Lifting Property (Theorem \ref{Path_Lifting_Property}), we have that $\widetilde{\overline{f \gamma_y} \cdot f \beta_y} = \widetilde{\overline{f \gamma_y}} \cdot \widetilde{f \beta_y}$. Hence,  $\widetilde{f \beta_y}(p_{\beta_y}^1) = \widetilde{f \gamma_y}(p_{\gamma_y}^1)$, which implies that $\widetilde{f}$ is well-defined.

The map $\widetilde{f}$ is also a graph homomorphism. Suppose the vertex $v$ is adjacent to $y$ in the graph $Y$. The map $\widetilde{f}$ is a graph homomorphism if there is an edge $\widetilde{f}(y) \widetilde{f}(v) \in E(\tX)$. Define $\alpha: I_\infty \to Y$ by
$$\alpha(i) = \begin{cases}
\gamma_y(i) & \text{for} \;\;\; i \leq p_{\gamma_y}^1\\
v &  \text{for} \;\;\; i \geq p_{\gamma_y}^1.
\end{cases}$$
It can be verified that $\alpha$ is a stable graph homomorphism with $p_\alpha^1 = p_{\gamma_y}^1 + 1$. Therefore, $\alpha(n_\alpha^1) = y_0$ and $\alpha(p_\alpha^1) = v$, which implies that $\widetilde{f}(v) = \widetilde{f \alpha}(p_\alpha^1)$. Since $\alpha(p_\alpha^1 - 1) = \alpha(p_{\gamma_y}^1) = \gamma_y(p_{\gamma_y}^1) = y$ and $f$ is a graph homomorphism, it follows that $f(\alpha(p_\alpha^1 - 1)) = f(\gamma_y(p_{\gamma_y}^1))$. Thus $\widetilde{f \gamma_y}(p_{\gamma_y}^1) = \widetilde{f \alpha}(p_\alpha^1 - 1)$, which implies that $\widetilde{f \gamma_y}(p_{\gamma_y}^1) \widetilde{f \alpha}(p_{\alpha}^1) \in E(\tX).$ In other words, $\widetilde{f}(y)\widetilde{f}(v)$ is an edge in $\tX$.

Finally, the map $\widetilde{f}$ is a lift of $f$. Since $\widetilde{f \gamma_y}: I_\infty \to \tX$ is a lift of $f \gamma_y$ and $\rho \circ \widetilde{f \gamma_y} = f \gamma_y$, it follows that
$$\rho \circ \widetilde{f}(y) = \rho(\widetilde{f \gamma_y}(p_{\gamma_y}^1)) = f(\gamma_y(p_{\gamma_y}^1)) = f(y)$$
for all $y \in V(Y)$. 
\end{proof}

In \cite{StromHomotopy}, the continuous maps which satisfy the Homotopy Lifting Property are one part of a Quillen model structure on topological spaces. Thus the Homotopy Lifting Property (Theorem \ref{Homotopy_Lifting_Property}) and the Lifting Criterion (Theorems \ref{Lifting_Criterion}) are a promising sign that A-homotopy theory might have a model structure associated to it. This would mean that A-homotopy theory is a homotopy theory for graphs in the Quillen sense.

\section{Application to $C_{k}$ for $k \geq 5$}
In classical homotopy theory, the first way we see covering spaces and the lifting properties used is in proving that the fundamental group of the circle is isomorphic to $\Z$. In \cite{BarceloFoundations}, Barcelo et al.  create a space from a graph by attaching 2-cells to the 3-cycles and 4-cycles of a graph. In Proposition 5.12 of \cite{BarceloFoundations}, Barcelo et al. prove that computing the classical fundamental group of this space is equivalent to computing the A-homotopy fundamental group of the original graph. This implies that the A-homotopy fundamental group of the $k$-cycle for $k \geq 5$ is isomorphic to $\Z$, since the spaces associated to these cycles will be homotopy equivalent to the circle. Similarly, the A-homotopy fundamental groups of $C_3$ and $C_4$ are isomorphic to $\{0\}$, since the spaces associated to these cycles is homotopy equivalent to the disk. In this section, we provide an alternate way to compute the A-homotopy theory fundamental groups of cycles $C_k$ for $k \geq 5$ analogous to the computation of the fundamental group of the circle in classical topology. This method involves covering graphs and the lifting properties. 

Computing the fundamental group of $C_k$ for $k \geq 5$ in this way gives one verification that the covering graphs and lifting properties work as expected in A-homotopy theory. In classical homotopy theory, covering maps are special cases of Hurewicz fibrations and an important part of a model structure on topological spaces. In future, we plan to investigate whether or not the covering graph maps are fibrations in a model structure on graphs. 




\begin{nota}\label{Covering_Graph_of_Ck_Nota}
Suppose that $k \geq 5$ and $C_k$ be a $k$-cycle with vertices labelled $[0], [1], \cdots, [k - 1]$, and let $\rho_k: I_\infty \to C_k$ be the graph homomorphism defined by $\rho_k(i) = [i \mod k]$ for all $i \in \Z$.
\end{nota}

Note that $\rho_k$ does not stabilize in either direction. 
For each $i \in V(I_\infty)$, we have that the neighborhood $N[i] = \{i-1, i, i+1\}$. The relative graph homomorphism $\rho_k|_{N[i]}$ is bijective for all $i \in \Z$. Thus $\rho_k$ is a local isomorphism, and the pair $(I_\infty, \rho_k)$ forms a covering graph of $C_k$.



\begin{ex}[Path Lift]\label{Path_Lift_Ex}
Let $\alpha \in B_1(C_k, [0])$, and let the pair $(I_\infty, \rho_k)$ be as in Notation \ref{Covering_Graph_of_Ck_Nota} By the Path Lifting Property (Theorem \ref{Path_Lifting_Property}), there is a unique lift $\widetilde{\alpha}: I_\infty \to I_\infty$ defined by $\widetilde{\alpha}(i) = 0$ for all $i \leq n_\alpha^1$, and recursively by $\widetilde{\alpha}(i) = (\rho_k|_{N_{\widetilde{\alpha}(i-1)}})^{-1}(\alpha(i))$ for all $i > n_\alpha^1$. 
\[\begin{tikzcd}
& & (I_\infty, \widetilde{x}) \arrow[dd, "p_k"]\\
& & \\
(I_\infty, 0) \arrow[rr, "\alpha"'] \arrow[rruu, "\widetilde{\alpha}"] & & (C_k, [0])
\end{tikzcd}\]
Since $N[\widetilde{\alpha}(i-1)] = \{\widetilde{\alpha}(i - 1) - 1, \widetilde{\alpha}(i - 1), \widetilde{\alpha}(i - 1) + 1\}$ and $\rho_k \circ \widetilde{\alpha} = \alpha$, 

\begin{eqnarray*}
\rho_k(\widetilde{\alpha}(i-1) - 1) & = & \alpha(i-1) - [1], \\
\rho_k(\widetilde{\alpha}(i-1)) & = & \alpha(i-1), \\
\rho_k(\widetilde{\alpha}(i-1) + 1) & = & \alpha(i-1) + [1].
\end{eqnarray*}
Thus
\begin{eqnarray*}
(\rho_k|_{N_{\widetilde{\alpha}(i-1)}})^{-1}(\alpha(i-1) - [1]) & = & \widetilde{\alpha}(i-1) - 1, \\
(\rho_k|_{N_{\widetilde{\alpha}(i-1)}})^{-1}(\alpha(i-1)) & = & \widetilde{\alpha}(i-1), \\
(\rho_k|_{N_{\widetilde{\alpha}(i-1)}})^{-1}(\alpha(i-1) + [1]) & = & \widetilde{\alpha}(i-1) + 1,
\end{eqnarray*}
and the lift $\widetilde{\alpha}: I_\infty \to I_\infty$ is defined for all $i \leq n_\alpha^1$ by $\widetilde{\alpha}(i) = 0$ and for all $i > n_\alpha^1$ recursively by
$$\widetilde{\alpha}(i) = 
\begin{cases}
\widetilde{\alpha}(i-1) + 1 & \text{if} \;\;\; \alpha(i) = \alpha(i-1) + [1], \\
\widetilde{\alpha}(i-1) & \text{if} \;\;\; \alpha(i) = \alpha(i-1), \\
\widetilde{\alpha}(i-1) - 1 & \text{if} \;\;\; \alpha(i) = \alpha(i-1) - [1].
\end{cases}$$
\end{ex}

We now propose an equivalence class of $A_1(C_k, [0])$ for each $n \in \Z$. We will later show that this is all of the equivalence classes of $A_1(C_k, [0])$.

\begin{df}\label{Gamma_n_Def}
Let $k \geq 5$ and the map $\gamma_n: I_\infty \to C_k$ be defined for each $n \geq 0$ by
$$\gamma_n(i) = 
\begin{cases}
[0] & \text{for} \;\;\; i \leq 0, \\
[i \mod k] & \text{for} \;\;\; 0 \leq i \leq kn, \\
[0] & \text{for} \;\;\; i \geq kn,
\end{cases}$$
and for each $n \leq 0$ by
$$\gamma_n(i) = 
\begin{cases}
[0] & \text{for} \;\;\; i \leq 0, \\
[(-i) \mod k] & \text{for} \;\;\; 0 \leq i \leq -kn, \\
[0] & \text{for} \;\;\; i \geq -kn,
\end{cases}$$
\end{df}

When $n = 0$, $\gamma_n$ is the constant map at $[0]$. For $n > 0$, the graph homomorphism $\gamma_n$ starts at $[0]$ and wraps around $C_k$ in a clockwise direction $n$ times. When $n < 0$, the graph homomorphism $\gamma_n$ starts at $[0]$ and wraps around $C_k$ in a counterclockwise direction $n$ times. By Example \ref{Path_Lift_Ex}, the lift $\widetilde{\gamma}_n: I_\infty \to I_\infty$ of $\gamma_n$ starting at $0$ is defined by
$$\widetilde{\gamma}_n(i) = 
\begin{cases}
0 & \text{for} \;\;\; i \leq 0, \\
i & \text{for} \;\;\; 0 \leq i \leq kn,\\
kn & \text{for} \;\;\; i \geq kn,
\end{cases}$$
if $n \geq 0$ and
$$\widetilde{\gamma_n}(i) = 
\begin{cases}
0 & \text{for} \;\;\; i \leq 0, \\
-i & \text{for} \;\;\; 0 \leq i \leq  -kn, \\
kn & \text{for} \;\;\; i \geq -kn,
\end{cases}$$
if $n \leq 0$. The Shifting Lemma (\ref{Shifting_Lemma}) can be used to relate the representatives $\gamma_n$ to each other. Recall that following Definition \ref{Fundamental_Group_Def}, the graph homomorphism $\overline{\gamma_n}$ is defined by $\overline{\gamma_n}(i) = \gamma_n(-i)$ for all $i \in \Z$. 

\begin{lemma}\label{Gamma_n_Inverse_Lemma}
Let $k \geq 5$ and $\gamma_n, \gamma_{-n} \in B_1(C_k, [0])$ be as defined in Definition \ref{Gamma_n_Def} for $n \in \Z$. Then $\gamma_{-n} \sim \overline{\gamma_n}$, whose equivalence class is the inverse of $[\gamma_n]$.
\end{lemma}

\begin{proof}
Suppose $n \geq 0$. By Definition \ref{Gamma_n_Def}, 
$$\overline{\gamma}_n(i) = \gamma_n(-i) = \begin{cases}
[0] & \text{for} \;\;\; i \leq -kn, \\
[(-i) \mod k] & \text{for} \;\;\; -kn \leq i \leq 0, \\
[0] & \text{for} \;\;\; i \geq 0.
\end{cases}$$
By the construction of $\gamma_n$,
\begin{eqnarray*}
\gamma_{-n}(i + kn)
& = & \begin{cases}
[0] & \text{for} \;\;\; i + kn \leq 0, \\
[(-i - kn) \mod k] & \text{for} \;\;\; 0 \leq i + kn \leq kn, \\
[0] & \text{for} \;\;\; i + kn \geq kn,
\end{cases} \\
& = & \begin{cases}
[0] & \text{for} \;\;\; i \leq - kn, \\
[(-i) \mod k] & \text{for} \;\;\; - kn \leq i  \leq 0, \\
[0] & \text{for} \;\;\; i \geq 0,
\end{cases} \\
& = & \overline{\gamma_n}(i),
\end{eqnarray*}
for all $i \in \Z$. Therefore, the graph homomorphism $\overline{\gamma_n}$ is equal to $\gamma_{-n}$ shifted by $-kn$. By the Shifting Lemma (\ref{Shifting_Lemma}), this implies that $\gamma_{-n} \sim \overline{\gamma_n}$ for $n \geq 0$. Similarly, for $n \leq 0$.
\end{proof}

Before proceeding to the proof that $A_1(C_k, [0]) \cong \Z$, we need to know which maps are A-homotopic to $\widetilde{\gamma_n}$. The analogous question in classical homotopy theory is solved using linear path homotopies. Unfortunately, the discrete nature of our paths requires a more complicated graph homotopy.

\begin{df}\label{Increasing_Decreasing_Def}
Let $\widetilde{f}: I_\infty \to I_\infty$ be a stable graph homomorphism. For $i \in \Z$, the value $\widetilde{f}(i)$ is \textit{increasing} if $\widetilde{f}(i) < \widetilde{f}(i + 1)$ and is \textit{decreasing} if $\widetilde{f}(i) > \widetilde{f}(i + 1)$ and is \textit{constant} if $\widetilde{f}(i) = \widetilde{f}(i + 1)$.
\end{df}

\begin{lemma}\label{Gamma_n_Lift_Equivalence_Class_Lemma}
Let $k \leq 3$. If $\widetilde{f}: I_\infty \to I_\infty$ is a stable graph homomorphism with $\widetilde{f}(n_{\widetilde{f}}^1) = 0$ and $\widetilde{f}(p_{\widetilde{f}}^1) = kn$, then $[\widetilde{f}] = [\widetilde{\gamma}_n]$, where $\widetilde{\gamma}_n$ is a lift of $\gamma_n: I_\infty \to C_k$.
\end{lemma}

\begin{proof}
Let $\widetilde{f}: I_\infty \to I_\infty$ be a stable graph homomorphism with $\widetilde{f}(n_{\widetilde{f}}^1) = 0$ and $\widetilde{f}(p_{\widetilde{f}}^1) = kn$ with $n \in \Z$. Although the path $\widetilde{f}$ starts at 0 and ends at $kn$, $\widetilde{f}$ may increase, decrease, or remain constant from the vertex $n_{\widetilde{f}}^1$ to the vertex $p_{\widetilde{f}}^1$. In contrast, for $n \geq 0$, $\widetilde{\gamma_n}$ increases constantly from starting at 0 to ending at $kn$, and for $n \leq 0$, $\widetilde{\gamma_n}$ decreases constantly from starting at 0 to ending at $kn$. We show that $\widetilde{f}$ is homotopic to $\widetilde{\gamma}_n$ by first showing that $\widetilde{f}$ is homotopic to a path that has no negative increasing values and no positive decreasing values.

Define $H: I_\infty \square I_\infty \to I_\infty$ for all $j \leq 0$ by $H(i, j) = \widetilde{f}(i)$, and recursively for all $j > 0$ by
$$H(i,j) = 
\begin{cases}
H(i, j-1) - 1 & \text{if} \;\;\; 0 \leq H(i + 1, j - 1) < H(i, j - 1), \\
H(i, j-1) & \text{if} \;\;\; 0 \leq H(i, j - 1) \leq H(i + 1, j - 1), \\
H(i, j-1) + 1 & \text{if} \;\;\; H(i, j - 1) < H(i + 1, j - 1) \leq 0, \\
H(i, j-1) & \text{if} \;\;\; H(i + 1, j - 1) \leq H(i, j - 1) \leq 0.
\end{cases}$$
First, we must confirm that these are all of the cases. Define $H_j: I_\infty \to I_\infty$ by $H_j(i) = H(i, j)$ for all $i, j \in \Z$. The first case is if $H_{j-1}(i)$ is a positive decreasing value. The second case is if $H_{j-1}(i)$ is a non-negative increasing or constant value. The third case is if $H_{j-1}(i)$ is a negative increasing value. The fourth case is if $H_{j-1}(i)$ is a non-positive decreasing or constant value. These are all possible cases. Note that the second and fourth cases overlap when $H_{j-1}(i) = 0$ and is a constant value. The map $H$ is well-defined, however, since $H(i, j) = H(i, j - 1)$ in both cases. It is routine to verify that $H$ is a graph homomorphism.

We now show that $H$ is stable. Since $H$ is a graph homomorphism, the restriction $H_j$ is a graph homomorphism. Since $\widetilde{f}$ is a stable graph homomorphism, the difference between $n_{\widetilde{f}}^1$ and $p_{\widetilde{f}}^1$ is finite. Thus there are a finite number of $m \in \Z$ with $n_{\widetilde{f}}^1 \leq m \leq p_{\widetilde{f}}^1$.
\begin{itemize}
    \item[(1)] Suppose $H_j(m) = 0$. By definition of $H$, either $0 = H_j(m) \leq H_j(m+1)$ and $H_{j+1}(m) = H_j(m)$, or $H_j(m + 1) \leq H_j(m) = 0$ and $H_{j+1}(m) = H_j(m)$. Thus if $H_j(m) = 0$, then $H_{j+1}(m) = 0$. This also implies that if $H_0(m) = \widetilde{f}(m) > 0$, then $H_j(m) \geq 0$ for all $j \geq 0$, and if $H_0(m) = \widetilde{f}(m) < 0$, then $H_j(m) \leq 0$ for all $j \geq 0$.
    \item[(2)] Suppose $H_j(m) > 0$. By definition of $H$, either $0 \leq H_j(m+1) < H_j(m)$ and $H_{j+1}(m) = H_j(m)-1$, or $0 \leq H_j (m) \leq H_j(m+1)$ and $H_{j+1}(m) = H_j(m)$. Thus $H_j(m)$ is constant or decreasing.
    \item[(3)] Suppose $H_j(m) < 0$. By definition of $H$, either $H_j(m) < H_j(m+1) \leq 0$ and $H_{j+1}(m) = H_j(m) + 1$, or $H_j(m+1) \leq H_j(m) \leq 0$ and $H_{j+1}(m) = H_j(m)$. Thus $H_j(m)$ is constant or increasing.
\end{itemize}

Observe that if there exists $j \geq 0$ such that $H_{j+1}(i) = H_j(i)$ for all $i \in \Z$, then $H$ stabilizes in the positive direction on the $2^{nd}$-axis, that is, the integer $p_H^2$ exists. For each $j \geq 0$, $H$ does not stabilize at $j$ in the positive direction in the $2^{nd}$-axis if and only if there exists some $m \in \Z$ with $n_{\widetilde{f}}^1 \leq m \leq p_{\widetilde{f}}^1$ such that $H_j(m) \neq H_{j+1}(m)$. We now count how many times it is possible for $H_j(m) \neq H_{j+1}(m)$ for $j \geq 0$ and $n_{\widetilde{f}}^1 \leq m \leq p_{\widetilde{f}}^1$. There are at most $p_{\widetilde{f}}^1 - n_{\widetilde{f}}^1$ choices for $m \in \Z$ with $n_{\widetilde{f}}^1 \leq m \leq p_{\widetilde{f}}^1$. By parts (1)-(3), for each such $m$, there are at most $|\widetilde{f}(m)|$ times that $H_j(m) \neq H_{j+1}(m)$. This implies that $H$ must stabilize in the positive direction on the $2^{nd}$-axis at a maximum of $j= \sum_m |\widetilde{f}(m)| < \infty$. Therefore, the integer $p_H^2$ exists.

It is now routine to show that $H$ is a path homotopy from $\widetilde{f}$ to $H_{p_H^2}$ as in Definition \ref{Path_Graph_Homotopy_Def}. By definition of $H$, the graph homomorphism $H_{p_H^2}$ has no positive decreasing value and no negative increasing values. Since $H_{p_H^2}$ starts at 0 as well, if $n \geq 0$, no negative increasing values implies that $H_{p_H^2}$ has no negative values at all. Thus no positive decreasing values implies that $H_{p_H^2}$ is constant or increasing from 0 to $kn$. Similarly, if $n \leq 0$, no positive decreasing values implies that $H_{p_H^2}$ has no positive values at all. Thus no negative increasing values implies that $H_{p_H^2}$ is constant or decreasing from 0 to $kn$. 
Thus by the Padding Lemma (\ref{Padding_Lemma}), $H_{p_H^2} \sim \widetilde{\gamma_n}$. Therefore, $\widetilde{f} \sim \widetilde{\gamma_n}$ for all $n \in \Z$. 
\end{proof}

We conclude this section by computing the fundamental group of the $k$-cycle for $k \geq 5$ using the lifting properties. The original computation of this group is a consequence of \cite[Proposition 5.12]{BarceloFoundations}. We offer a new proof here.

\begin{prop}\cite[Proposition 5.12]{BarceloFoundations}\label{Fundamental_Group_Ck_Thm}
The fundamental group of $C_k$ for $k \geq 5$ is $(A_1(C_k, [0]), \cdot) \cong (\Z, +)$.
\end{prop}

\begin{proof}
Define $\varphi: \Z \to A_1(C_k, [0])$ by $n \mapsto [\gamma_n]$, the homotopy class of the stable graph homomorphism $\gamma_n: I_\infty \to C_k$ defined in Definition \ref{Gamma_n_Def}. The function $\varphi$ is a group homomorphism by the following cases.
\begin{itemize}
    \item Case 1: Suppose $n, m \geq 0$. Then $\gamma_{n+m}(i) = (\gamma_n \cdot \gamma_m)(i-km)$ for all $i \in \Z$, and $\gamma_{n+m} \sim \gamma_n \cdot \gamma_m$ by the Shifting Lemma (\ref{Shifting_Lemma}).
    \item Case 2: Suppose $n, m < 0$. Then $\gamma_{n+m}(i) = (\gamma_n \cdot \gamma_m)(i + km)$ for all $i \in \Z$, and $\gamma_{n+m} \sim \gamma_n \cdot \gamma_m$ by the Shifting Lemma (\ref{Shifting_Lemma}).
    \item Case 3: Suppose $n \geq 0, m < 0$. By Lemma \ref{Gamma_n_Inverse_Lemma}, $\overline{\gamma_n} \sim \gamma_{-n}$ and $\overline{\gamma_m} \sim \gamma_{-m}$. By Case 1, if $n + m \geq 0$, then $\gamma_n = \gamma_{n+m-m} \sim \gamma_{n+m} \cdot \gamma_{-m}$. By Case 2, if $n + m < 0$, then $\gamma_m = \gamma_{-n+n+m} \sim \gamma_{-n} \cdot \gamma_{n+m}$. Thus
    $$\gamma_n \cdot \gamma_m \sim \gamma_n \cdot \overline{\gamma_{-m}} \sim \gamma_{n+m} \cdot \gamma_{-m} \cdot \overline{\gamma_{-m}} \sim \gamma_{n+m} \;\;\;\;\; \text{if} \;\;\; n+m \geq 0,$$
    and 
    $$\gamma_n \cdot \gamma_m \sim \overline{\gamma_{-n}} \cdot \gamma_m \sim \overline{\gamma_{-n}} \cdot \gamma_{-n} \cdot \gamma_{n+m} \sim \gamma_{n+m} \;\;\;\;\; \text{if} \;\;\; n+m \leq 0.$$
    \item Case 4: Suppose that $n < 0, m \geq 0$. Similar to Case 3.
\end{itemize}
Therefore, $\varphi(n+m) = [\gamma_{n+m}] = [\gamma_n \cdot \gamma_m] = [\gamma_n] \cdot [\gamma_m] = \varphi(n) \cdot \varphi(m)$ for all $n, m \in \Z$.

The function $\varphi$ is also surjective. If $[f] \in A_1(C_k, [0])$, then $f$ is a stable graph homomorphism with $f(n_f^1) = f(p_f^1) = [0]$. Hence, by the Path Lifting Property (Theorem \ref{Path_Lifting_Property}), there exists a unique lift $\widetilde{f}: I_\infty \to I_\infty$ with $\widetilde{f}(n_f^1) = 0$ and $f = \rho_k \circ \widetilde{f}$. Since $f(p_f^1) = [0]$, it follows that $\rho_k(\widetilde{f}(p_f^1)) = [0]$, so 
$$\widetilde{f}(p_f^1) \mod k = 0.$$
Thus there exists $n \in \Z$ such that $\widetilde{f}(p_f^1) = kn$. Hence, by the Lemma \ref{Gamma_n_Lift_Equivalence_Class_Lemma}, we have that $\widetilde{f} \sim \widetilde{\gamma_n}$, which implies that there exists a path homotopy $H: I_\infty^2 \to I_\infty$ from $\widetilde{f}$ to $\widetilde{\gamma_n}$. To conclude, verify that $\rho_k \circ H$ is a path homotopy from $f$ to $\gamma_n$ as in Definition \ref{Path_Graph_Homotopy_Def}. 

Finally, $\varphi$ is injective. If $\varphi(n) = \varphi(m)$, then $[\gamma_n] = [\varphi_m]$, which implies that $\gamma_n \sim  \gamma_m$. Therefore, there exists a path homotopy $H: I_\infty^2 \to C_k$ from $\gamma_n$ to $\gamma_m$. By the Homotopy Lifting Property (Theorem \ref{Homotopy_Lifting_Property}), there is a path homotopy $\widetilde{H}: I_\infty^2 \to I_\infty$ from $\widetilde{\gamma_n}$ to $\widetilde{\gamma_m}$. Thus $\widetilde{\gamma_n} \sim \widetilde{\gamma_m}$, which implies that $\widetilde{\gamma_n}(p_{\varphi_n}^1) = \widetilde{\gamma_m}(p_{\gamma_m}^1)$. Therefore, $kn = km$, and it follows that $n = m$.

Thus $\varphi$ is an isomorphism, and $A_1(C_k, [0]) \cong \Z$.
\end{proof}


\begin{thebibliography}{}
\bibitem{Atkin1} R.Atkin, \textit{An Algebra for Patterns on a Complex, I}, Internat. J. Man–Machine Studies \textbf{6} (1974), pp. 285–307.

\bibitem{Atkin2} R. Atkin, \textit{An Algebra for Patterns on a Complex, II}, Internat. J. Man–Machine Studies \textbf{8} (1976), pp. 448–483.

\bibitem{BabsonHomotopy} E. Babson, H. Barcelo, M. de Longueville, and R. Laubenbacher, \textit{Homotopy Theory of Graphs}. J. Algebraic Combinatorics \textbf{24} (2006), no. 1, pp. 31–44. https://arxiv.org/pdf/math/0403146.pdf

\bibitem{BarceloFoundations} H. Barcelo, X. Kramer, R. Laubenbacher, C. Weaver, \textit{Foundations of a Connectivity Theory for Simplicial Complexes}. Advances in Applied Mathematics \textbf{26} (2001), no. 2, pp. 97–128.

\bibitem{BarceloPerspectives}  H. Barcelo, R. Laubenbacher, \textit{Perspectives on A-Homotopy Theory and its Applications}. Discrete Mathematics \textbf{298} (2005), no. 1-3, pp. 39–61. 

\bibitem{GodsilAlgebraic} C. Godsil, G. F. Royle, \textit{Algebraic Graph Theory}, Graduate Texts in Mathematics \textbf{207}. Springer Verlag, New York, (2001).

\bibitem{HardemanMscThesis} R. Hardeman, \textit{ Computing A-Homotopy Groups of Graphs Using Coverings and Lifting Properties}, Master's thesis, University of Calgary. http://hdl.handle.net/1880/108868, (2018).

\bibitem{HatcherAlgebraic} A. Hatcher, \textit{Algebraic Topology}, Cambridge University Press, New York, 2001. 

\bibitem{HellGraphs} P. Hell, J. Nesetril, \textit{Graphs and Homomorphisms}, Oxford University Press, New York, 2004.

\bibitem{MunkresTopology} J. R. Munkres, \textit{Topology}, 2nd ed., Prentice Hall, Inc., Upper Saddle River, NJ, 2004.

\bibitem{StromHomotopy} A. Strom, \textit{The Homotopy Category is a Homotopy Category}, Archiv der Mathematik \textbf{23} (1973), no. 23, 435-441

\bibitem{WestGraph} D. West, \textit{ Introduction to Graph Theory}, 2nd ed., Prentice Hall, Inc., Upper Saddle River, NJ,  2001.
\end{thebibliography}
\end{document}